# AN APPROXIMATE SAMPLING FORMULA UNDER GENETIC HITCHHIKING

By Alison Etheridge,[1] Peter Pfaffelhuber[2] and
Anton Wakolbinger[3]

*University of Oxford, Ludwig-Maximilian University Munich
and Goethe-University Frankfurt*

For a genetic locus carrying a strongly beneficial allele which has just fixed in a large population, we study the ancestry at a linked neutral locus. During this "selective sweep" the linkage between the two loci is broken up by recombination and the ancestry at the neutral locus is modeled by a structured coalescent in a random background. For large selection coefficients $\alpha$ and under an appropriate scaling of the recombination rate, we derive a sampling formula with an order of accuracy of $\mathcal{O}((\log \alpha)^{-2})$ in probability. In particular we see that, with this order of accuracy, in a sample of fixed size there are at most two nonsingleton families of individuals which are identical by descent at the neutral locus from the beginning of the sweep. This refines a formula going back to the work of Maynard Smith and Haigh, and complements recent work of Schweinsberg and Durrett on selective sweeps in the Moran model.

**1. Introduction.** Assume that part of a large population of size $2N$ carries, on some fixed genetic locus (henceforth referred to as the *selective locus*), an allele with a certain selective advantage. If the population reproduction is described by a classical Fisher–Wright model or, more generally, a Cannings model with individual offspring variance $\sigma^2$ per generation, and time is measured in units of $2N$ generations, the evolution of the fraction

Received March 2005; revised September 2005.

[1]Supported by EPSRC Advanced Fellowship GR/A90923.

[2]Travel support from DFG, Bilateral Research Group FOR 498.

[3]Supported in part by a Senior Research Fellowship at Erwin Schroedinger Institute, Vienna.

carrying the advantageous allele is approximately described by the Fisher–Wright stochastic differential equation (SDE)

$$(1.1) \qquad dP = \sqrt{\sigma^2 P(1-P)}\, dW + \alpha P(1-P)\, dt,$$

where $W$ is a standard Wiener process and $s = \alpha/2N$ is the selective advantage of the gene per individual per generation [6, 8].

Assume at a certain time a sample of size $n$ is drawn from the subpopulation carrying the advantageous allele. Conditioned on the path $P$, the ancestral tree of the sample *at the selective locus* is described by Kingman's coalescent with pair coalescence rate $\sigma^2/P$ (see, e.g., [14] and Remark 4.6 below).

Now consider a *neutral locus* in the neighborhood of the selective one, with a recombination probability $r$ per individual per generation between the two loci. From generation to generation, there is a small probability $r$ per individual that the gene at the neutral locus is broken apart from its selective partner and *recombined* with another one, randomly chosen from the population. In the diffusion limit considered here, this translates into the *recombination rate* $\rho$. Depending on $P$, only a fraction of these recombination events will be *effective* in changing the status of the selective partner from "advantageous" to "nonadvantageous" or vice versa. Given $P$, the genealogy of the sample at the neutral locus can thus be modeled by a *structured coalescent of the neutral lineages in background $P$* as in [2]: Backward in time, a neutral lineage currently linked to the advantageous allele recombines with a nonadvantageous one at rate $\rho(1-P)$, and a neutral lineage currently linked to a nonadvantageous gene recombines with the advantageous one at rate $\rho P$, where $\rho = 2Nr$. Moreover, two neutral lineages currently both linked to the advantageous allele coalesce at rate $\sigma^2/P$, and two neutral lineages currently both linked to a nonadvantageous allele coalesce at rate $\sigma^2/(1-P)$.

Two individuals sampled at time $t > 0$ are said to be *identical by descent at the neutral locus from time* 0 if their neutral ancestral lineages coalesce between times 0 and $t$. This defines an *ancestral sample partition at the neutral locus* at time $t$ from time 0.

We are interested in a situation in which at a certain time (say time 0) a single copy of the advantageous gene enters into the population and eventually fixates. For large $N$, the time evolution of the size of the subpopulation carrying the advantageous allele can be thought of as a random path $X$ governed by an $h$-transform of (1.1), entering from $x = 0$ and conditioned to hit $x = 1$.

The parameters $\alpha$, $\rho$ and $\sigma^2$, the random path $X$, its fixation time $T$ and the structured $n$-coalescent in background $X$ from time $T$ back to time 0 are the principal ingredients in the first part of our analysis. Our central object of interest is the ancestral sample partition at the neutral locus at time $T$ from time 0.

For simplicity we put $\sigma^2 = 2$. This not only simplifies some formulae, but also allows a direct comparison with the results of Schweinsberg and Durrett [21], who considered the finite population analog in which the population evolves according to a Moran model, leading to $\sigma^2 = 2$ in the diffusion approximation.

We focus on large coefficients $\alpha$ and refer to $X$ as the *random path in a selective sweep*. The expected duration of the sweep is approximately $(2 \log \alpha)/\alpha$ (see Lemma 3.1). Heuristically, the sweep can be divided into several *phases*. Phases that must certainly be considered are: the time intervals which $X$ takes:

- to reach a small level $\varepsilon$ (phase 1);
- to climb from $\varepsilon$ to $1 - \varepsilon$ (phase 2);
- to fixate in 1 after exceeding the level $1 - \varepsilon$ (phase 3).

Whereas the expected durations of phases 1 and 3 both are $\asymp \log \alpha / \alpha$, that of phase 2 is only $\asymp 1/\alpha$. The analysis of hitchhiking has in the past often concentrated on the second phase [15, 22]. For large population size and large selective advantage, the frequency path of the beneficial allele in phase 2 can be approximately described by a deterministic logistic growth curve (see, e.g., [17]). However, this approximation is only good for fixed $\varepsilon > 0$. In [15] the frequency path under a selective sweep is described by a logistic growth curve that model phase 2 with $\varepsilon = 5/\alpha$ (recall that $\alpha$ stands for $2Ns$), whereas in [22] $\varepsilon = 1/(2N)$ is considered. In both cases, $\varepsilon$ decreases with population size. As Barton pointed out in [1], the logistic model fails to include the randomness of the frequency path at the very beginning of the selective sweep. Consequently he further subdivided phase 1 so as to study the onset of the sweep in more detail.

No matter how the phases of a sweep are chosen, we have seen that the first phase as given above takes $\asymp \log \alpha / \alpha$. So, to see a nontrivial number of recombination events along a single lineage between $t = 0$ and $t = T$, the recombination rate $\rho$ should be on the order of $\alpha / \log \alpha$. Henceforth, we therefore assume

$$(1.2) \qquad \rho = \gamma \frac{\alpha}{\log \alpha}, \qquad 0 < \gamma < \infty.$$

With this recombination rate, it will turn out that, asymptotically as $\alpha \to \infty$, effectively no recombinations happen in phase 2, since this phase is so short; neither do they occur in phase 3, since then $1 - X$ is so small. Consequently, the probability that a single ancestral lineage is not hit by a recombination is approximately given by

$$(1.3) \qquad p = e^{-\gamma}.$$

Since for large $\alpha$ the subpopulation carrying the advantageous allele is expanding quickly near time $t = 0$, a first approximation to the sample genealogy at the selective locus is given by a star-shaped tree, that is, $n$ lineages all coalescing at $t = 0$. Hence, ignoring possible back-recombinations (which can be shown to have small probability; see Proposition 3.4), a first approximation to the sampling distribution at the neutral locus is given by a Binomial$(n, p)$ number of individuals stemming from the founder of the sweep; the rest of the individuals are recombinants that all have different neutral ancestors at the beginning of the sweep (see Remark 2.6 and cf. [21], Theorem 1.1). This approximate sampling formula goes back to the pioneering work of Maynard Smith and Haigh [18], who also coined the term *hitchhiking*: the allele which the founder of the sweep carried at the neutral locus gets a lift into a substantial part of the population (and the sample) in the course of the sweep. Apart from the hitchhikers, there are also a few free riders who jump on the sweep when it is already under way and make it as singletons into the sample.

We will see that the sampling formula just described is accurate with probability $1 - \mathcal{O}(\frac{1}{\log \alpha})$. In a spirit similar to [21], but with a somewhat different strategy, we will improve the order of accuracy up to $\mathcal{O}(\frac{1}{(\log \alpha)^2})$. A common technical theme is an approximation of the genealogy of the advantageous allele in the *early phase* of the sweep by a supercritical branching process [1, 21], an idea which can be traced back to Fisher [9] (see [8], page 27ff). It is this early phase of the sweep which is relevant for the deviations from the Maynard Smith and Haigh sampling distribution. This higher-order approximation allows for the possibility of the occurrence of recombination events during this early phase that affect our sample and thus, potentially, lead to nonsingleton families associated to a neutral type other than the original hitchhiker.

Our main result is the derivation of a sampling formula for the ancestral partition at the end of the sweep which is accurate up to an error $\mathcal{O}(1/(\log \alpha)^2)$ in probability.

## 2. Model and main result.

### 2.1. *The model.*

DEFINITION 2.1 (Random path of a selective sweep). For $\alpha > 0$, let $X = (X_t)_{0 \le t \le T}$ be a random path in $[0, 1]$ following the SDE

$$(2.1) \quad dX = \sqrt{2X(1 - X)}\, dW + \alpha X(1 - X) \coth\left(\frac{\alpha}{2}X\right) dt, \qquad 0 < t,$$

and entering from 0 at time $t = 0$. Here, $W$ is a standard Wiener process and $T$ denotes the time when $X$ hits 1 for the first time.

Note that 0 is an entrance boundary for (2.1) and that $X$ given by (2.1) arises as an $h$-transform of the solution of (1.1). Indeed, since the latter has generator

$$Gf(x) = x(1-x)f''(x) + \alpha x(1-x)f'(x)$$

and the $G$ harmonic function $h : [0,1] \to [0,1]$ with boundary conditions $h(0) = 0$, $h(1) = 1$ is

$$h(x) = \frac{1 - e^{-\alpha x}}{1 - e^{-\alpha}},$$

the $h$-transformed generator is

$$G^h f(x) = \frac{1}{h(x)} G(hf)(x)$$

$$= x(1-x)f''(x) + \alpha x(1-x) \coth\left(\frac{\alpha}{2}x\right)f'(x);$$

see also [11], page 245. As described in the Introduction, $X$ models the evolution of the size of the subpopulation that consists of all those individuals which carry an advantageous allele (called $B$) on a *selective locus*. At time $t = 0$, a single mutant that carries the (then novel) allele $B$ enters into the population, corresponding to $X_0 = 0$, and $X$ is conditioned to eventually hit 1, which happens at the random time $T$. We will refer to $X$ given by (2.1) as the *random path of the sweep* (or *random sweep path* for short).

DEFINITION 2.2 (*$n$-coalescent in background $X$*). Let $n \geq 2$. Given a random sweep path $X$, we construct a tree $\mathcal{T}_n = \mathcal{T}_n^X$ with $n$ leaves as follows. Attach the $n$ leaves of $\mathcal{T}_n^X$ at time $t = T$ and work backward in time $t$ from $t = T$ to $t = 0$ (or equivalently, forward in time $\beta = T - t$ from $\beta = 0$ to $\beta = T$) with pair coalescence rate $2/X_t$ [i.e., starting from the tree top ($\beta = 0$), the *current number of lineages* in $\mathcal{T}_n^X$ decreases at rate $\frac{2}{X_t}\binom{n}{2}$]. Finally, attach the root of $\mathcal{T}_n^X$ at time $t = 0$.

Note that $\mathcal{T}_n^X$ corresponds to a *time-changed Kingman coalescent* (see [12]). This time change transforms the one single lineage of infinite length, which can be thought to follow the ultimate coalescence in Kingman's coalescent, into a single lineage ending at $\beta = T$. We will refer to the tree $\mathcal{T}_n^X$ as the *$n$-coalescent in background $X$*; it describes the *genealogy at the selective locus* of a sample of size $n$ taken from the population at the instant of completion of the sweep.

DEFINITION 2.3 (*Coalescing Markov chains in background $X$*). Let $\rho > 0$. Given a random sweep path $X$, let $(\xi_\beta)_{0 \leq \beta \leq T}$ be a $\{B, b\}$-valued Markov chain with time inhomogeneous jump rates

$$(2.2) \qquad \begin{aligned} (1 - X_{T-\beta})\rho &= (1 - X_t)\rho &\quad \text{from } B \text{ to } b, \\ X_{T-\beta}\rho &= X_t\rho &\quad \text{from } b \text{ to } B. \end{aligned}$$

The process $\xi$ describes to which type at the selective locus (either $B$ or $b$) an ancestral lineage of the neutral locus is currently linked in its journey back into the past, indexed by the backward time $\beta$. Recall that each *neutral gene* (i.e., each gene at the neutral locus) is linked at any time to a *selective gene* (i.e., at the selective locus), the latter being either of type $B$ or of type $b$, and that for the neutral lineages which are currently in state $B$, only a fraction $1 - X_t$ of the recombination events is effective in taking them into state $b$.

DEFINITION 2.4 (Structured $n$-coalescent in background $X$; see [2]). Given $X$, consider $n$ independent copies of the Markov chain $\xi$ (see Definition 2.3), all of them starting in state $B$ at time $\beta = 0$. Let any two $\xi$-walkers who are currently (say at time $\beta = T - t$) in state $B$ coalesce at rate $2/X_{T-\beta} = 2/X_t$ and let any two $\xi$-walkers in state $b$ coalesce at rate $2/(1 - X_t)$. The resulting (exchangeable) system $\Xi_n$ of coalescing Markov chains will be called the *structured $n$-coalescent in background $X$*.

We now define a *labeled partition* $\mathcal{P}^{\Xi_n}$ of $\{1, \ldots, n\}$ induced by the structured coalescent $\Xi_n$.

DEFINITION 2.5 (Ancestral sample partition at the neutral locus). In the situation of Definition 2.4, we will say that $i$ and $j$, $1 \le i < j \le n$, *belong to the same family* if the $\xi$-chains numbered $i$ and $j$ in $\Xi_n$ coalesce before time $\beta = T$.

1. A family will be labeled *nonrecombinant* if none of the ancestral lineages of the family back to time $\beta = T$ ever left state $B$. (Thus the neutral ancestor of a nonrecombinant family at time $t = 0$ is necessarily linked to the single selective founder of the sweep.)
2. A family will be labeled *early recombinant* if none of its ancestral lineages ever left state $B$ before the first (looking backward from $t = T$) coalescence in the sample genealogy happened, but if nonetheless the family's ancestor at time $t = 0$ is in state $b$.
3. A family will be labeled *late recombinant* if at least one of its ancestral lineages left state $B$ before the first (looking backward from $T$) coalescence in the sample genealogy happened *and* if the family's ancestor at time $t = 0$ is in state $b$.

4. In all other cases (e.g., if two lineages on their way back first leave $B$, then coalesce and return to $B$ afterward), the family will be labeled *exceptional*.

The *labeled partition* resulting in this way will be called $\mathcal{P}^{\Xi_n}$.

For large selection coefficients and moderately large recombination rates [see (1.2)] it turns out that, up to an error in probability of $\mathcal{O}(\log \alpha)^{-2}$, all late recombinant families are singletons, there is no more than one early recombinant family and there are no exceptional families. In fact, the probability that there is an early recombinant family at all is of the order $(\log \alpha)^{-1}$. Given there is an early recombinant family, however, its size may well be a substantial fraction of $n$. Our main result (Theorem 1) clarifies the approximate distribution of the number of late recombinants and of the size of the early recombinant family.

2.2. *Main result.* Recall that $\mathcal{P}^{\Xi_n}$ (introduced in Definition 2.5) describes the ancestral partition of an $n$-sample drawn from the population at the time of completion of the sweep, where the partition is induced by identity by descent at the neutral locus at the beginning of the sweep.

**THEOREM 1** (Approximate distribution of the ancestral sample partition). *Fix a sample size $n$. For a selection coefficient $\alpha \gg 1$ and a recombination rate $\rho$ obeying (1.2) for fixed $\gamma$, the random partition $\mathcal{P}^{\Xi_n}$ introduced in Definition 2.5 consists, with probability $1 - \mathcal{O}((\log \alpha)^{-2})$, of the following parts:*

- *$L$ late recombinant singletons.*
- *One family of early recombinants of size $E$.*
- *One nonrecombinant family of size $n - L - E$.*

*More precisely, a random labeled partition of $\{1, \ldots, n\}$, whose distribution approximates that of $\mathcal{P}^{\Xi_n}$ up to a variation distance of order $\mathcal{O}((\log \alpha)^{-2})$, is given by random numbers $L$ and $E$ constructed as follows:*

*Let $F$ be an $\mathbb{N}$-valued random variable with*

$$(2.3) \qquad \mathbf{P}[F \leq i] = \frac{(i - (n-1)) \cdots (i-1)}{(i + (n-1)) \cdots (i+1)},$$

*and, given $F = f$, let $L$ be a binomial random variable with $n$ trials and success probability $1 - p_f$, where*

$$(2.4) \qquad p_f = \exp\left(-\frac{\gamma}{\log \alpha} \sum_{i=f}^{\lfloor \alpha \rfloor} \frac{1}{i}\right).$$

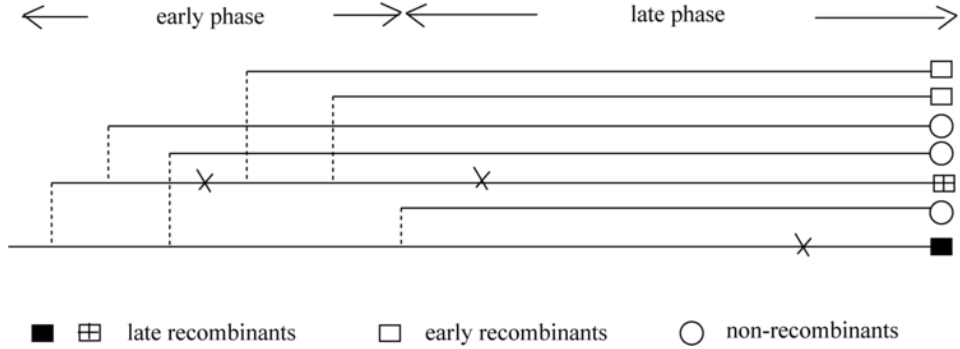

Fig. 1.

*Independently of all this, let $S$ be a $\{0, 1, \ldots, n\}$-valued random variable with*

$$(2.5) \qquad \mathbf{P}[S = s] = \begin{cases} \dfrac{\gamma n}{\log \alpha} \displaystyle\sum_{i=2}^{n-1} \dfrac{1}{i}, & s = 1, \\[2ex] \dfrac{\gamma n}{\log \alpha} \dfrac{1}{s(s-1)}, & 2 \le s \le n-1, \\[2ex] \dfrac{\gamma n}{\log \alpha} \dfrac{1}{n-1}, & s = n. \end{cases}$$

*Given $S = s$ and $L = l$, the random variable $E$ is hypergeometric, choosing $n - l$ out of $n = s + (n - s)$, that is,*

$$(2.6) \qquad \mathbf{P}[E = e] = \frac{\binom{s}{e}\binom{n-s}{n-l-e}}{\binom{n}{n-l}}.$$

Sections 3 and 4 will be devoted to the proof of Theorem 1. Figure 1 explains the concepts which appear in our theorem and points to the strategy of the proof as explained in Section 3.2. In the figure, the sample size is $n = 7$, and the $\times$'s indicate effective recombination events that occur along the lineages. The early phase ends when the number of lines in the sample tree has increased from six to seven. At the end of the early phase there is one family of size $S = 3$ of early recombinants. One member of this family is then kicked out by a late recombination. In the sample there are $L = 2$ late recombinant singletons, one early recombinant family of size $E = 2$ and one nonrecombinant family of size 3.

REMARK 2.6. (a) From (2.5) we see that

$$\mathbf{P}[S \ge 2] = \frac{\gamma n}{\log \alpha}$$

and

$$\mathbf{P}[S > 0] = \frac{\gamma n}{\log \alpha} \sum_{i=1}^{n-1} \frac{1}{i}.$$

In particular, this shows that the probability that there are early recombinants at all is $\mathcal{O}(\frac{1}{\log \alpha})$.

(b) Barton [1] reported simulations in which several nonsingleton recombinant families arise. This is not, in fact, incompatible with our theorem since the constant in the $\mathcal{O}$ in the error estimate of Theorem 1 depends on $(\gamma n)^2$. [See, e.g., Section 3.4, where for each pair in our $n$-sample we encounter an error $C \frac{\gamma^2}{(\log \alpha)^2}$.] In, for example, the simulation described on page 130 of [1], in which eight early recombinant families are seen, this factor $\gamma n$ is $\approx 120$, while $\log \alpha \approx 13$, which explains the occurrence of several nonsingleton recombinant families.

As a corollary to Theorem 1 we obtain an approximate sampling formula under the model for genetic hitchhiking. This means that we can now derive the probability of having $l$ late recombinants (which produce singletons), $e$ early recombinants, which form a family of size $e$, and $n - l - e$ lineages that go back to the founder of the sweep and also form a family on their own.

COROLLARY 2.7 (Approximate sampling formula). *Under the assumptions of Theorem* 1 *the common distribution of the number of early recombinants* $E$ *and the number of late recombinants* $L$ *is, with probability* $1 - \mathcal{O}((\log \alpha)^{-2})$, *given by*

$$\mathbf{P}[E = e, L = l]$$

$$(2.7) \qquad = \mathbf{E}[p_F^{n-l}(1 - p_F)^l] \cdot \begin{cases} \dfrac{n\gamma}{\log \alpha} \dfrac{(n-1)\binom{n-2}{e-2}\mathbb{1}\{l + e = n\} + \binom{n-1}{l}}{e(e-1)}, \\ \hfill e \geq 2, \\[2mm] \dfrac{n\gamma}{\log \alpha}\Bigg(\mathbb{1}\{l + 1 = n\} \\ \qquad + \binom{n-1}{l}\sum_{i=2}^{n-1}\dfrac{1}{i} + \sum_{s=2}^{n}\dfrac{\binom{n-s}{l-s+1}}{s-1}\Bigg), \\ \hfill e = 1, \\[2mm] \binom{n}{l}\Bigg(1 - \dfrac{n\gamma}{\log \alpha}\Bigg(1 - \dfrac{l}{n}\sum_{i=2}^{n-1}\dfrac{1}{i}\Bigg)\Bigg) \\ \qquad + \dfrac{n\gamma}{\log \alpha}\Bigg(\dfrac{1}{n}\mathbb{1}\{l = n\} \\ \qquad\qquad + \sum_{s=2}^{n}\binom{n-s}{n-l}\dfrac{1}{s(s-1)}\Bigg), \\ \hfill e = 0. \end{cases}$$

*Necessarily given $L = l$ and $E = e$, the number of lineages going back to the founder of the sweep is $n - l - e$.*

This corollary will be proved in Section 4.5.

2.3. *Comparison with Schweinsberg and Durrett's work.* Our research has been substantially inspired by recent work of Schweinsberg and Durrett [4, 21]. Let us point out briefly how the main results of [21] and of the present paper complement each other.

Schweinsberg and Durrett [21] considered a two-locus Moran model with population size $2N$, selective advantage $s$ of the advantageous allele and individual recombination probability $r = \mathcal{O}(1/\log N)$. Their main result is (in our terminology) about the approximate distribution of the ancestral distribution of an $n$-sample at the neutral locus as $N \to \infty$. In preparing their Theorem 1.2, Schweinsberg and Durrett [21] specified (in terms of a stick-breaking scheme made up by a sequence of Beta variables) a random paintbox with parameters $L = \lfloor 2Ns \rfloor$ and $r/s$. They denoted the (labeled) distribution of an $n$-sample drawn from the paintbox (where the class belonging to the first draw is tagged) by $Q_{r/s,L}$. The assertion of their Theorem 1.2 is then that $Q_{r/s,L}$ approximates the ancestral sample distribution at the neutral locus with probability $1 - \mathcal{O}(1/(\log N)^2)$. Notably, $s$ remains fixed, that is, does *not* scale with $N$.

A priori, this *strong selection limit* does not lend itself to a diffusion approximation. However, interestingly enough, certain aspects do: in particular, the ones studied in the present context. More precisely, our results show that the approximate distribution of the ancestral sample partition in the strong selection limit of [21] arises also in a two-stage way, first passing to the diffusion limit and then letting the selection coefficient tend to infinity. This is made precise in the following proposition, which will be proved at the end of Section 4.4.

PROPOSITION 2.8. *Let $Q_{r/s,L}$ be as in [21], Theorem 1.2. Then, with the choice*

$$(2.8) \qquad \alpha = 2Ns, \qquad \rho = 2Nr, \qquad \gamma = \frac{r}{s}\log\alpha,$$

*the distribution specified in Theorem 1 and further described in Corollary 2.7 approximates $Q_{r/s,L}$ up to an error of $\mathcal{O}(1/(\log N)^2)$.*

Hence our results (Theorem 1 and Corollary 2.7) also give an approximate sampling formula for the random partition that appears in Theorem 1.2 of [21] and enters as an input to the coalescent with simultaneous multiple mergers described in [5]. In particular, our results reveal that this random partition has more than one nonsingleton class with a probability of only $\mathcal{O}(1/(\log N)^2)$, a result which is less explicit in the proofs of [21].

Let us emphasize once again that in [21] the error is controlled in a specific "large but finite population" model, whereas in our approach the error is controled after having performed a diffusion approximation. Proposition 2.8 together with Theorem 1.2 of [21] reveals that our diffusion approximation has the order of approximation $\mathcal{O}(1/(\log N)^2)$ in the strong selection limit of the Moran model. This might be seen as one more indication for the strength and robustness of the diffusion approximation in the context of population genetics.

*Numerical results.* One can now still ask how large are the constants which are lurking behind the $\mathcal{O}$'s. To shed some light on this and to see how well our approximations perform, let us present some numerics. We compare the approximation of Theorem 1 with numerical examples given in [21]. The examples deal with samples of size $n = 1$ and $n = 2$. We distinguish the number and types of ancestors of the sample at the beginning of the sweep. For a single individual the probability that the ancestor is of type $b$ (an event called pinb in [21]) can be approximated by

$$\text{pinb} \approx \mathbf{P}[L = 1],$$

as there is no early phase in our theorem in this case. For a sample of size 2, either there are two ancestors and both have type $b$ (denoted p2inb), there are two ancestors, one of type $B$ and one of type $b$ (denoted p1B1b), or there is one ancestor with either a $b$ allele (denoted p2cinb) or a $B$ allele (which happens in all other cases). Using Theorem 1 we approximate

$$\text{p2inb} \approx \mathbf{P}[L = 2 \text{ or } S = 2, L = 1],$$

$$\text{p2cinb} \approx \mathbf{P}[L = 0, S = 2],$$

$$\text{p1B1b} \approx \mathbf{P}[L = 1, S = 0].$$

In [21] simulations were performed for three models: (i) in a Moran model, (ii) in a model where the frequency of the $B$ allele follows a deterministic logistic growth curve, and (iii) for the approximate result obtained in [21], Theorem 1.2. Results for a more extensive range of parameters can be found in [4]. In Table 1, we have added the approximations of our Theorem 1 to those of [21]. In all cases we find that the approximation given by our Theorem 1 performs comparably to the approximation of Schweinsberg and Durrett. Both approximations are significantly better than the logistic model.

**3. Outline of the proof of Theorem 1.** We start by calculating the expected duration of the sweep.

TABLE 1
*Numerical results comparing Theorem 1 with Theorem 2 from [21] and a logistic model.*
*The numbers in brackets are relative errors with respect to the Moran model*

|            | **pinb**      | **p2inb**     | **p2cinb**     | **p1B1b**      |
|------------|---------------|---------------|----------------|----------------|
|            | $n = 10^4$    | $s = 0.1$     | $r = 0.001064$ |                |
| Moran      | 0.08203       | 0.00620       | 0.01826        | 0.11513        |
| Logistic   | 0.09983(21%)  | 0.00845(36%)  | 0.03365(84%)   | 0.11544(0.3%)  |
| DS, Thm. 2 | 0.08235(0.4%) | 0.00627(1.1%) | 0.01765(−3.4%) | 0.11687(1.5%)  |
| Thm. 1     | 0.08249(0.6%) | 0.00659(6.3%) | 0.01867(2.2%)  | 0.11515(0.0%)  |
|            | $n = 10^4$    | $s = 0.1$     | $r = 0.005158$ |                |
| Moran      | 0.33656       | 0.10567       | 0.05488        | 0.35201        |
| Logistic   | 0.39936(18%)  | 0.13814(31%)  | 0.09599(75%)   | 0.32646(−7.3%) |
| DS, Thm. 2 | 0.34065(1.2%) | 0.10911(3.2%) | 0.05100(−7.1%) | 0.36112(2.6%)  |
| Thm. 1     | 0.32973(−2.0%)| 0.10857(2.7%) | 0.05662(3.2%)  | 0.34157(−3.0%) |

3.1. *The duration of the sweep.* Let $T_\delta$ be the time at which $X$ reaches the level $\delta$ for the first time.

LEMMA 3.1. *For all fixed $\varepsilon \in (0,1)$, as $\alpha \to \infty$,*

$$(3.1) \qquad \mathbf{E}[T_\varepsilon] = \frac{\log \alpha}{\alpha} + \mathcal{O}\left(\frac{1}{\alpha}\right), \qquad \mathbf{E}[T - T_{1-\varepsilon}] = \frac{\log \alpha}{\alpha} + \mathcal{O}\left(\frac{1}{\alpha}\right),$$

$$(3.2) \quad \mathbf{E}[T_{1-\varepsilon} - T_\varepsilon] = \mathcal{O}\left(\frac{1}{\alpha}\right),$$

$$(3.3) \qquad \mathbf{Var}[T] = \mathcal{O}\left(\frac{1}{\alpha^2}\right).$$

This lemma will be proved in Section 4. Notice in particular that $\mathbf{E}T = 2\log \alpha / \alpha + \mathcal{O}(\alpha^{-1})$. Thus, to see a nontrivial number of recombination events along a single line between $t = 0$ and $t = T$ ($\beta = T$ and $\beta = 0$), the recombination rate $\rho$ should be on the order of $\alpha / \log \alpha$. Henceforth, we therefore assume that $\rho$ obeys equation (1.2).

3.2. *Three approximation steps.* In a *first approximation step* we will show that all events (both coalescences and recombinations) that happen along the lineages of the structured coalescent while dwelling in state $b$ have a negligible effect on the sampling distribution. Thus once a lineage has recombined away from state $B$, we can assume that it experiences no further recombination or coalescence events in the remaining time to $\beta = T$. This motivates us to couple the structured $n$-coalescent (at the neutral locus) with the $n$-coalescent (at the selective locus), and to study the ancestral partition at the neutral locus by marking effective recombination events that happen along the selective lineages.

DEFINITION 3.2.   For a given sweep path $X$, consider the coalescent $\mathcal{T}_n$ in background $X$ together with a Poisson process with intensity measure $(1 - X_t)\rho \, dt$ along the lineages of $\mathcal{T}_n$. Say that two leaves of $\mathcal{T}_n$ belong to the same family if and only if the path in $\mathcal{T}_n$ which connects them is not hit by a mark. Call a mark *early* if it occurs between time 0 and the time when $\mathcal{T}_n$ increases from $n - 1$ to $n$ and call it *late* otherwise. Label a family as *early-marked* (*late-marked*) if it traces back to an early (late) mark; otherwise label it as *nonmarked*. In this way we arrive at what we call the *labeled partition* $\mathcal{P}^{\mathcal{T}_n}$.

Note that the nonlabeled version of $\mathcal{P}^{\mathcal{T}_n}$ arises from the marked tree $\mathcal{T}_n$ in the same way as the sample partition in the infinite-alleles model emerges from the marked coalescent. Also note that late-marked families in $\mathcal{P}^{\mathcal{T}_n}$ are necessarily singletons. It will turn out (see Corollary 3.5) that $\mathcal{P}^{\mathcal{T}_n}$ approximates $\mathcal{P}^{\Xi_n}$ with probability $1 - \mathcal{O}((\log \alpha)^{-2})$.

The *second approximation step* consists of replacing $\mathcal{P}^{\mathcal{T}_n}$ by a labeled partition $\mathcal{P}^{\mathcal{Y}_n}$ generated by a marked Yule tree.

DEFINITION 3.3.   Let $\mathcal{Y}$ be an (infinite) Yule tree with branching rate $\alpha$ and let $\mathcal{Y}_n$ be the random tree which arises by sampling $n$ lineages from $\mathcal{Y}$ (which come down from infinity). Up to the time when the number of lines extant in $\mathcal{Y}$ reaches the number $\lfloor \alpha \rfloor$, mark the lines of $\mathcal{Y}_n$ by a Poisson process with homogeneous intensity $\rho = \gamma \alpha / \log \alpha$. Families, early and late marks, early-marked families and so forth are specified in complete analogy to Definition 3.2. The resulting labeled partition of $\{1, \ldots, n\}$ will be denoted by $\mathcal{P}^{\mathcal{Y}_n}$.

Again, we will show (see Proposition 3.6) that $\mathcal{P}^{\mathcal{Y}_n}$ approximates $\mathcal{P}^{\mathcal{T}_n}$ with probability $1 - \mathcal{O}((\log \alpha)^{-2})$.

The *third approximation step* exploits the fact that the probability for more than one early mark is $\mathcal{O}((\log \alpha)^{-2})$. The random variable $F$ specified in Theorem 1 stands for the number of lines extant in the full tree $\mathcal{Y}$ at the time of the most recent coalescence in the sample tree $\mathcal{Y}_n$. The number $p_f$ is the approximate probability that, given $F = f$, a single lineage is not hit by a late mark (or equivalently, does not experience a late recombination); note that this probability is larger than the probability $p$ given by (1.3). The number of late-marked families (corresponding to late recombinants) is approximated by a mixed Binomial random variable $L$ with random success probability $1 - p_F$. In the dominating case that $\mathcal{Y}_n$ in its early phase (when it has less than $n$ lines) is hit by at most one mark, the random variable $S$ approximates the size of the early-marked family which arises if we "cut off" $\mathcal{Y}_n$ at the time of its most recent coalescence (i.e., when its number of lines increases from $n - 1$ to $n$). This size is thinned out by late marks; in other words, the final size of the early-marked family arises as a hypergeometric random variable, by randomly distributing the $n - L$ lineages which have not been knocked out by late marks onto the $S + (n - S)$ potential ancestors at the most recent coalescence time of $\mathcal{Y}_n$.

3.3. *First step*: *From the structured to a marked coalescent.* The key result for the first approximation step is the following:

PROPOSITION 3.4.  (i) *The probability that a neutral ancestral lineage of our sample recombines out of $B$ and then recombines back into $B$ and* (ii) *the probability that a pair of neutral ancestral lineages coalesces in $b$ are both* $\mathcal{O}(\frac{1}{(\log \alpha)^2})$.

The previous proposition allows us [within the accuracy of $\mathcal{O}(\frac{1}{(\log \alpha)^2})$] to dispense with the structured coalescent and work instead with the *marked coalescent* in background $X$ as described in Definition 3.2. Indeed, the following statement is immediate from Proposition 3.4.

COROLLARY 3.5.   *The variation distance between the distributions of* $\mathcal{P}^{\Xi_n}$ *and* $\mathcal{P}^{\mathcal{T}_n}$ *is* $\mathcal{O}(\frac{1}{(\log \alpha)^2})$.

3.4. *Second step*: *From the marked coalescent to a marked Yule tree.* A key tool will be a time transformation which takes the random sweep path into a (stopped) supercritical Feller diffusion. Because the early phase of the sweep is the most relevant one and because a Fisher–Wright diffusion entering from zero looks similar to a Feller diffusion entering from zero as long as both are close to zero, we will be able to control the error in the sample partition that results from replacing $X$ by the path of a Feller diffusion. Under this time transformation the mark (recombination) rate will become a constant. By exchangeability, we can (without loss of generality) sample only from individuals in our Feller diffusion with infinite line of descent (provided of course that there are enough such) and it is well known that such "immortal particles" form a Yule tree with branching rate $\alpha$ [7]. This means that we sample from a Poisson number of individuals with parameter $\alpha$, but we shall see that it suffices to consider the Yule tree stopped when it has precisely $\lfloor \alpha \rfloor$ extant individuals.

PROPOSITION 3.6.   *Let* $\mathcal{Y}$, $\mathcal{Y}_n$ *and* $\mathcal{P}^{\mathcal{Y}_n}$ *be as in Definition 3.3. The variation distance between the distributions of* $\mathcal{P}^{\mathcal{T}_n}$ *and* $\mathcal{P}^{\mathcal{Y}_n}$ *is* $\mathcal{O}(\frac{1}{(\log \alpha)^2})$.

3.5. *Third step*: *Approximating sample partitions in marked Yule trees.*
Because of Corollary 3.5 and Proposition 3.6, the proof of Theorem 1 will
be complete if we can show that the representation given there applies,
within the accuracy of $\mathcal{O}(\frac{1}{(\log \alpha)^2})$, to the random labeled partition $\mathcal{P}^{\mathcal{Y}_n}$.
Thus, the remaining part of the proof takes place in the world of marked Yule
processes, where matters are greatly simplified and many exact calculations
are possible.

Let $I = I(t)$ be the number of lines of $\mathcal{Y}$ extant at time $t$ and let $K_i$ be
the number of lines extant in $\mathcal{Y}_n$ while $I = i$. The process $K = (K_i)$ will play
a major role in our analysis below. Viewing the index $i$ as *time*, referred to
below as *Yule time*, we will see that $K$ is a Markov chain.

We denote by $M_i$ the number of marks that hit $\mathcal{Y}_n$ while $\mathcal{Y}$ has $i$ lines.
Since the latter period is exponentially distributed with parameter $i\alpha$ and
marks appear along lines according to a Poisson process with rate $\gamma\alpha/\log \alpha$,
we arrive at the following observation:

REMARK 3.7.   Given $K_i = k$, $M_i$ is distributed as $G - 1$, where $G$ has a
geometric distribution with parameter

$$\frac{i\alpha}{i\alpha + k\gamma\alpha/\log \alpha} = \frac{1}{1 + (k/i)\gamma/\log \alpha}.$$

Consequently, the conditional expectation of $M_i$ given $K_i = k$ is $\frac{k}{i}\frac{\gamma}{\log \alpha}$.

We will distinguish two phases of the process $\mathcal{Y}$. The *early phase* will
consist of all Yule times $i$ when $K_i < n$ and the *late phase* will consist of the
Yule times $i$ with $K_i = n$.

We define

$$F := \min\{i : K_i = n\},$$

that is, $F$ is the number of lines in the full tree $\mathcal{Y}$ when the number of lines
in the sample tree $\mathcal{Y}_n$ reaches its final size $n$. This is when the late phase
begins.

PROPOSITION 3.8.   *The distribution of $F$ is given by* (2.3).

By analogy with Definition 3.2, we call those marks which hit $\mathcal{Y}_n$ in the
early (late) phase the *early* (*late*) *marks*.

The labeled partition $\mathcal{P}^{\mathcal{Y}_n}$ introduced in Definition 3.3 is generated by
early and late marks. Let us treat the early marks first. We will find in
Proposition 3.9 that up to our desired accuracy, there is at most one early
mark.

We write

$$M = \sum_{i\,:\,i < F} M_i$$

for the number of early marks. Let us denote by $S^{\mathcal{Y}}$ the number of leaves in $\mathcal{Y}_n$ whose ancestral lineage is hit by an early mark. On the event $\{M = 1\}$, that is, in the case of a single early mark, the leaves of $\mathcal{T}_n$ are partitioned into two classes, one (of size $S^{\mathcal{Y}}$) whose ancestry is hit by this single early mark and one whose ancestry is not hit by an early mark.

The next proposition gives an approximation for the joint distribution of $(M, S^{\mathcal{Y}})$.

PROPOSITION 3.9.   *Up to an error of $\mathcal{O}(\frac{1}{(\log \alpha)^2})$,*

$$(3.4) \qquad \mathbf{P}[M = 1, S^{\mathcal{Y}} = s] = \begin{cases} \dfrac{n\gamma}{\log \alpha} \sum_{k=2}^{n-1} \dfrac{1}{k}, & s = 1, \\[2ex] \dfrac{n\gamma}{\log \alpha} \dfrac{1}{s(s-1)}, & 2 \leq s \leq n-1, \\[2ex] \dfrac{n\gamma}{\log \alpha} \dfrac{1}{n-1}, & s = n. \end{cases}$$

*Furthermore,*

$$(3.5) \qquad \mathbf{P}[M \geq 2] = \mathcal{O}\left(\frac{1}{(\log \alpha)^2}\right).$$

For fixed $f \leq \lfloor \alpha \rfloor$, the probability that a randomly chosen line is not hit by a mark between $i = f$ and $i = \lfloor \alpha \rfloor$, is (cf. Remark 3.7)

$$\begin{aligned} \prod_{i=f}^{\lfloor \alpha \rfloor} \frac{1}{1 + (1/i)\gamma/\log \alpha} &= \exp\left( \sum_{i=f}^{\lfloor \alpha \rfloor} \log\left( 1 - \frac{\gamma}{\log \alpha} \frac{1}{i + \gamma/\log \alpha} \right) \right) \\ (3.6) \qquad &= \exp\left( -\frac{\gamma}{\log \alpha} \sum_{i=f}^{\lfloor \alpha \rfloor} \frac{1}{i + \gamma/\log \alpha} \right) + \mathcal{O}\left( \frac{1}{(\log \alpha)^2} \right) \\ &= p_f + \mathcal{O}\left( \frac{1}{(\log \alpha)^2} \right), \end{aligned}$$

where $p_f$ was defined in (2.4). The last step follows from the Taylor expansion

$$\frac{1}{i + \gamma/\log \alpha} = \frac{1}{i} + \frac{1}{i^2}\mathcal{O}\left(\frac{1}{\log \alpha}\right).$$

Write $L^{\mathcal{Y}}$ for the number of lineages in $\mathcal{Y}_n$ that are hit by late marks. For distinct lineages of $\mathcal{Y}_n$, the events that they are hit by late marks are asymptotically independent, which allows us to approximate the distribution of $L^{\mathcal{Y}}$.

PROPOSITION 3.10. *The distribution of $L^{\mathcal{Y}}$ is approximately mixed Binomial. More precisely,*

$$
\begin{aligned}
(3.7) \quad \mathbf{P}[L^{\mathcal{Y}} = l] &= \sum_{f=n}^{\lfloor \alpha \rfloor} \binom{n}{l} p_f^{n-l} (1 - p_f)^l \mathbf{P}[F = f] + \mathcal{O}\left(\frac{1}{(\log \alpha)^2}\right) \\
&= \binom{n}{l} \mathbf{E}[p_F^l (1 - p_F)^{n-l}] + \mathcal{O}\left(\frac{1}{(\log \alpha)^2}\right), \qquad l = 0, \dots, n.
\end{aligned}
$$

Based on the previous two propositions, we will be able to show that up to our desired accuracy the random variables $S^{\mathcal{Y}}$ and $L^{\mathcal{Y}}$ can be treated as independent.

PROPOSITION 3.11. *The random variables $S^{\mathcal{Y}}$ and $L^{\mathcal{Y}}$ are approximately independent, that is,*

$$
\mathbf{P}[S^{\mathcal{Y}} = s, L^{\mathcal{Y}} = l] = \mathbf{P}[S^{\mathcal{Y}} = s] \cdot \mathbf{P}[L^{\mathcal{Y}} = l] + \mathcal{O}\left(\frac{1}{(\log \alpha)^2}\right).
$$

Given $M = 1, S^{\mathcal{Y}} = s$ and $L^{\mathcal{Y}} = l$, the size (call it $E^{\mathcal{Y}}$) of the (single) early-marked family in $\mathcal{P}^{\mathcal{Y}_n}$ (see Definition 3.3) is hypergeometric, choosing $n - l$ out of two classes, one of size $s$, the other of size $n - s$. Hence from Propositions 3.9, 3.10 and 3.11 the labeled partition $\mathcal{P}^{\mathcal{Y}_n}$ consists, with probability $1 - \mathcal{O}((\log \alpha)^{-2})$, of $L$ late-marked singletons, one early-marked family of size $E$ and one nonmarked family of size $n - L - E$, where the joint distribution of $(L, E)$ is specified in Theorem 1. Combining this with Proposition 3.6 and Corollary 3.5, the proof of Theorem 1 is complete.

## 4. Proofs.

4.1. *The duration of the sweep: Proof of Lemma 3.1.* We use standard theory about one-dimensional diffusions. (See, e.g., [8] and [16].) Primarily we need the Green's function $G(\cdot, \cdot)$ that corresponds to the solution of the SDE (2.1), this time with $X_0 = x \in [0, 1]$.

The Green's function $G(\cdot, \cdot)$ satisfies

$$
(4.1) \qquad \mathbf{E}_x\left[\int_0^T g(X_s)\, ds\right] = \int_0^1 G(x, \xi) g(\xi)\, d\xi,
$$

where $\mathbf{E}_x[\cdot]$ is the expectation with respect to the process $X$ started in $x$ (and $\mathbf{E}[\cdot] = \mathbf{E}_0[\cdot]$).

If $X$ is a solution of (2.1), then (see [16], Chapter 15, formula (9.8))

(4.2)
$$x \leq \xi : G(x,\xi) = \frac{(1-e^{-\alpha(1-\xi)})(1-e^{-\alpha\xi})}{\alpha\xi(1-\xi)(1-e^{-\alpha})},$$

$$x \geq \xi : G(x,\xi) = \frac{(e^{-\alpha x}-e^{-\alpha})(e^{\alpha\xi}-1)(1-e^{-\alpha\xi})}{\alpha\xi(1-\xi)(1-e^{-\alpha})(1-e^{-\alpha x})}.$$

Observe that $G(x,\xi)$ is decreasing in $x$.

PROOF OF (3.1) AND (3.2). In the proofs there will appear some constants $C, C'$ which might change from occurrence to occurrence.

Observe that

(4.3)          $$\mathbf{E}[T_\varepsilon] = \mathbf{E}[T] - \mathbf{E}_\varepsilon[T] = \int_0^\varepsilon G(0,\xi)\,d\xi - \int_0^\varepsilon G(\varepsilon,\xi)\,d\xi,$$

where we have used that $G(0,\xi) = G(\varepsilon,\xi)$ as long as $\xi \geq \varepsilon$. Since

$$\frac{1}{\xi(1-\xi)} = \frac{1}{\xi} + \frac{1}{1-\xi}$$

and using the symmetry in $G(0,\xi) = G(0,1-\xi)$, we see that for the first term in (4.3),

$$\int_0^\varepsilon G(0,\xi)\,d\xi - \frac{\log\alpha}{\alpha} = \frac{1}{\alpha}\left(\int_0^\varepsilon + \int_{1-\varepsilon}^1 \frac{(1-e^{-\alpha\xi})(1-e^{-\alpha(1-\xi)})}{(1-e^{-\alpha})\xi}\,d\xi - \log\alpha\right)$$

$$= \frac{1}{\alpha}\left(C + \int_0^{\alpha\varepsilon} \frac{(1-e^{-\xi})(1-e^{-\alpha+\xi})}{\xi}\,d\xi - \int_1^{\alpha\varepsilon}\frac{1}{\xi}\,d\xi\right)$$

$$= \frac{1}{\alpha}\left(C - \int_1^{\alpha\varepsilon}\frac{e^{-\xi}+e^{-\alpha+\xi}-e^{-\alpha}}{\xi}\,d\xi\right) = \mathcal{O}\left(\frac{1}{\alpha}\right).$$

For the second term, as $1 - e^{-\alpha\xi} \leq 1 - e^{-\alpha\varepsilon}$ for $\xi \leq \varepsilon$,

$$\int_0^\varepsilon G(\varepsilon,\xi)\,d\xi \leq \frac{e^{-\alpha\varepsilon}}{\alpha(1-\varepsilon)}\int_0^\varepsilon \frac{e^{\alpha\xi}-1}{\xi}\,d\xi = \frac{e^{-\alpha\varepsilon}}{\alpha(1-\varepsilon)}\int_0^{\alpha\varepsilon}\frac{e^\xi-1}{\xi}\,d\xi$$

$$\leq \frac{Ce^{-\alpha\varepsilon}}{\alpha}\left(1 + \int_1^{\alpha\varepsilon} e^\xi\,d\xi\right) = \mathcal{O}\left(\frac{1}{\alpha}\right).$$

A similar calculation leads to the second statement of (3.1), and from these two equalities and as $\mathbf{E}[T] = \mathbf{E}[T_{1/2}] + \mathbf{E}[T - T_{1/2}]$, also (3.2) follows. □

PROOF OF (3.3). To compute the variance of $T$ we use the following identity, which is a consequence of the Markov property (and can be checked

by induction on $k$):

$$(4.4) \quad \mathbf{E}_x\left[\int_0^T \int_{t_1}^T \cdots \int_{t_{k-1}}^T g_k(X_{t_k}) \cdots g_1(X_{t_1}) \, dt_k \cdots dt_1\right]$$

$$= \int_0^1 \cdots \int_0^1 G(x, x_1) \cdots G(x_{k-1}, x_k) g_1(x_1) \cdots g_k(x_k) \, dx_k \cdots dx_1.$$

From this we obtain

$$
\begin{aligned}
\mathbf{Var}[T] &= 2\int_0^1 \int_0^1 G(0, \xi) G(\xi, \eta) \, d\eta \, d\xi - 2\int_0^1 \int_\xi^1 G(0, \xi) G(0, \eta) \, d\eta \, d\xi \\[2mm]
&= 2\int_0^1 \int_0^\xi G(0, \xi) G(\xi, \eta) \, d\eta \, d\xi \\[2mm]
&= \frac{2}{\alpha^2(1 - e^{-\alpha})^2} \\[2mm]
&\quad \times \int_0^1 \int_0^\xi \frac{(1 - e^{-\alpha(1-\xi)})(e^{-\alpha\xi} - e^{-\alpha})}{\xi(1-\xi)} \frac{(e^{\alpha\eta} - 1)(1 - e^{-\alpha\eta})}{\eta(1-\eta)} \, d\eta \, d\xi \\[2mm]
&\overset{\xi \to 1-\xi}{\leq} \frac{2e^{-\alpha}}{\alpha^2(1 - e^{-\alpha})^2} \int_0^1 \int_0^{1-\xi} \frac{e^{\alpha\xi} - 1}{\xi(1-\xi)} \frac{e^{\alpha\eta} - 1}{\eta(1-\eta)} \, d\eta \, d\xi \\[2mm]
&= \frac{2}{\alpha^2(1 - e^{-\alpha})^2} \left(\left(e^{-\alpha/2} \int_0^{1/2} \frac{e^{\alpha\xi} - 1}{\xi(1-\xi)} \, d\xi\right)^2 \right. \\[2mm]
&\qquad\qquad \left. + 2\int_0^{1/2} \int_0^\xi \frac{e^{-\alpha\xi} - e^{-\alpha}}{\xi(1-\xi)} \frac{e^{\alpha\eta} - 1}{\eta(1-\eta)} \, d\eta \, d\xi\right)
\end{aligned}
$$

by a decomposition of the area $\{(\xi, \eta) : \eta \leq 1 - \xi\}$ in $\{(\xi, \eta) : \xi, \eta \leq 1/2\}$, $\{(\xi, \eta) : \xi \leq 1/2, 1/2 \leq \eta \leq 1 - \xi\}$ and $\{(\xi, \eta) : \xi \geq 1/2, \eta \leq 1 - \xi\}$, and the symmetry of the integrand. From this we see

$$\mathbf{Var}[T] \leq \frac{2}{\alpha^2}\left(C + 8\int_0^{1/2} \int_0^\xi \frac{e^{-\alpha\xi} - e^{-\alpha}}{\xi} \frac{e^{\alpha\eta} - 1}{\eta} \, d\eta \, d\xi\right)$$

and (3.3) follows as

$$
\begin{aligned}
\int_0^{1/2} &\int_0^\xi \frac{e^{-\alpha\xi} - e^{-\alpha}}{\xi} \frac{e^{\alpha\eta} - 1}{\eta} \, d\eta \, d\xi \\[2mm]
&\overset{\xi \to \alpha\xi}{\underset{\eta \to \alpha\eta}{=}} \int_0^{\alpha/2} \frac{e^{-\xi} - e^{-\alpha}}{\xi} \int_0^\xi \frac{e^\eta - 1}{\eta} \, d\eta \, d\xi \\[2mm]
&\leq \int_1^{\alpha/2} \frac{e^{-\xi} - e^{-\alpha}}{\xi} \left(\int_1^\xi \frac{e^\eta - 1}{\eta} \, d\eta + C\right) d\xi + C' \\[2mm]
&\leq C\int_1^\alpha \frac{1}{\xi^2} \, d\xi + C' = \mathcal{O}(1).
\end{aligned}
$$

Here we have used $\int_0^{\alpha/2} \int_0^\xi d\eta \, d\xi = \int_0^1 \int_0^\xi d\eta \, d\xi + \int_1^{\alpha/2}(\int_0^1 d\eta + \int_1^\xi d\eta) \, d\xi$ to obtain the penultimate and $\int_1^\xi d\eta = \int_1^{\xi/2} d\eta + \int_{\xi/2}^\xi d\eta$ to obtain the last inequality. $\quad\square$

### 4.2. *From the structured to a marked coalescent.*

PROOF OF PROPOSITION 3.4. (i) Given $(X_t)_{0 \le t \le T}$, we are looking for the probability that tracing backward from time $T$ to time $0$ a lineage escapes the sweep (which happens with rate $\rho(1 - X_t)$) and then recombines back into $B$ (which happens then with rate $\rho X_t$). The required probability then follows by integrating over path space and is given by

$$\mathbf{E}\left[\int_0^T \left(1 - \exp\left(-\int_0^t \rho X_s \, ds\right)\right) \rho(1 - X_t) \exp\left(-\int_t^T \rho(1 - X_s) \, ds\right) dt\right]$$

$$(4.5) \qquad \le \rho^2 \mathbf{E}\left[\int_0^T (1 - X_t) \int_0^t X_s \, ds \, dt\right]$$

$$\le \rho^2 \int_0^1 \int_0^\xi G(0, \xi) G(\xi, \eta) \, d\eta \, d\xi$$

$$+ \rho^2 \int_0^1 \int_\xi^1 G(0, \xi) G(0, \eta)(1 - \eta) \xi \, d\eta \, d\xi,$$

where we have used that $G(\xi, \eta) = G(0, \eta)$ for $\xi \le \eta$ and (4.4). The first term is $\frac{\rho^2}{2} \mathbf{Var}[T] = \mathcal{O}(\frac{1}{(\log \alpha)^2})$ by (3.1). The second term gives

$$\frac{\rho^2}{\alpha^2(1 - e^{-\alpha})^2} \int_0^1 \int_\xi^1 \frac{(1 - e^{-\alpha\xi})(1 - e^{-\alpha(1-\xi)})}{1 - \xi}$$

$$\times \frac{(1 - e^{-\alpha\eta})(1 - e^{-\alpha(1-\eta)})}{\eta} \, d\eta \, d\xi$$

$$\le \frac{\rho^2}{\alpha^2} \int_0^1 \int_0^\eta \frac{1}{(1 - \xi)\eta} \, d\xi \, d\eta = \frac{\rho^2}{\alpha^2} \int_0^1 \int_{1-\eta}^1 \frac{1}{\xi} \frac{1}{\eta} \, d\xi \, d\eta$$

$$= \frac{\rho^2}{\alpha^2} \int_0^1 -\frac{\log(1 - \eta)}{\eta} \, d\eta = \mathcal{O}\left(\frac{1}{(\log \alpha)^2}\right)$$

and we are done.

(ii) To prove the second assertion of the proposition, we split the event that two lineages coalesce in $b$ into two events. Recall that $T_\delta$ denotes the time when $X$ first hits $\delta$. Whenever two lineages coalesce in $b$, then either there must have been two recombination events between $T_{1/2}$ and $T$ or there must have been a coalescence in $b$ between $0$ and $T_{1/2}$. Both events only have small probabilities as we now show.

First consider the event that two recombinations occur in $[T_{1/2}, T]$. We see here as in (4.5) that the probability for this event is at most

$$\rho^2 \mathbf{E}_{1/2}\left[\int_0^T \int_0^t (1 - X_t)(1 - X_s)\, ds\, dt\right]$$

$$= \rho^2 \int_0^1 \int_0^1 G(\tfrac{1}{2}, \xi) G(\xi, \eta)(1-\eta)(1-\xi)\, d\eta\, d\xi$$

$$\leq \rho^2 \int_0^1 \int_0^\xi G(0, \xi) G(\xi, \eta)\, d\eta\, d\xi$$

$$+ \rho^2 \int_0^{1/2} \int_\xi^1 G(\tfrac{1}{2}, \xi) G(0, \eta)(1-\xi)(1-\eta)\, d\eta\, d\xi$$

$$+ \rho^2 \int_{1/2}^1 \int_{1/2}^1 G(0, \xi) G(0, \eta)(1-\xi)(1-\eta)\, d\eta\, d\xi,$$

where we have used that $G(x, \xi)$ is decreasing in $x$. The first term is $\frac{\rho^2}{2} \mathbf{Var}[T] = \mathcal{O}(\frac{1}{(\log \alpha)^2})$. The second is bounded by

$$\frac{\rho^2}{\alpha^2} \int_0^{1/2} \int_\xi^1 \frac{e^{-\alpha/2}(e^{\alpha \xi} - 1)}{\xi} \frac{1 - e^{-\alpha \eta}}{\eta}\, d\eta\, d\xi$$

$$= \frac{\gamma^2 e^{-\alpha/2}}{(\log \alpha)^2} \int_0^\alpha \int_0^{\alpha/2 \wedge \eta} \frac{e^\xi - 1}{\xi} \frac{1 - e^{-\eta}}{\eta}\, d\xi\, d\eta$$

$$\leq \frac{\gamma^2 e^{-\alpha/2}}{(\log \alpha)^2} \left(C + \int_1^\alpha \frac{e^{\eta \wedge \alpha/2}}{\eta}\, d\eta\right) = \mathcal{O}\left(\frac{1}{(\log \alpha)^2}\right),$$

where we have split the integral $\int_1^\alpha d\eta$ in $\int_1^{\alpha/2} d\eta + \int_{\alpha/2}^\alpha d\eta$ to obtain the final estimate. The third term is small, as it is the square of

$$\rho \int_{1/2}^1 G(0, \xi)(1-\xi)\, d\xi \leq \frac{\gamma}{\log \alpha} \int_{1/2}^1 \frac{1}{\xi}\, d\xi = \mathcal{O}\left(\frac{1}{\log \alpha}\right).$$

The second event we have to consider is coalescence in $b$ between time 0 and $T_{1/2}$. The probability for this event is, again as in (4.5), at most

$$\mathbf{E}\left[\int_0^{T_{1/2}} \frac{2}{1 - X_t}\, dt\right] \leq \mathbf{E}\left[\int_0^T \frac{2}{1 - X_t} \mathbf{1}\left(X_t \leq \frac{1}{2}\right) dt\right] = \int_0^{1/2} G(0, \xi) \frac{2}{1 - \xi}\, d\xi$$

$$\leq \frac{2}{\alpha} \int_0^{1/2} \frac{1 - e^{-\alpha \xi}}{\xi (1 - \xi)^2}\, d\xi \leq \frac{8}{\alpha} \int_0^{\alpha/2} \frac{1 - e^{-\xi}}{\xi}$$

$$\leq \frac{8}{\alpha}(C + \log \alpha) = \mathcal{O}\left(\frac{\log \alpha}{\alpha}\right).$$

So both events are improbable and Proposition 3.4(ii) is proved. $\quad \square$

4.3. *From the coalescent to the Yule tree.*

PROOF OF PROPOSITION 3.6.   We have to show that the error we make when using Yule trees instead of coalescent trees in a random background is small. This involves two approximation steps. The first is an approximation of the coalescent in the random background of a Wright–Fisher diffusion by a coalescent in an $\alpha$-supercritical Feller background. The second step is an approximation of the latter coalescent by Yule trees.

For the first approximation step we need to time-change our coalescent. This relies on the following proposition, whose proof consists of an application of [6], Chapter 6, Section 1.

PROPOSITION 4.1.   *Under the random time change $t \mapsto \tau$ given by $d\tau = (1 - X_t) \, dt$, the random path $X = (X_t)_{0 \le t \le T}$ is taken into a random path $Z = (Z_\tau)_{0 \le \tau \le \tilde{T}}$, which is an $\alpha$-supercritical Feller diffusion governed by*

$$dZ = \sqrt{2Z} \, dW + \alpha Z \, d\tau,$$

*starting in $Z_0 = 0$, conditioned on nonextinction and stopped at the time $\tilde{T}$ when it first hits 1.*

Under the time change $t \mapsto \tau$, the $n$-coalescent $\mathcal{T}_n$ described in Definition 2.2 is taken into the $n$-coalescent $\mathcal{C}_n$ whose pair coalescence rate conditioned on $Z$ is $\frac{2}{Z_\tau(1 - Z_\tau)} \, d\tau$. Under this time change, the marking rate $\rho(1 - X_t) \, dt$ becomes $\rho \, d\tau$. Let $\mathcal{P}^{\mathcal{C}_n}$ be the sample partition generated along $\mathcal{C}_n$ in the same way as $\mathcal{P}^{\mathcal{T}_n}$ was generated along $\mathcal{T}_n$, but now with the uniform marking rate $\rho \, d\tau$. Note that $\mathcal{P}^{\mathcal{C}_n}$ and $\mathcal{P}^{\mathcal{T}_n}$ have the same distribution.

Let us denote by $\mathcal{D}_n$ the $n$-coalescent whose pair coalescence rate conditioned on $Z$ is $2/Z_\tau \, d\tau$, $\tilde{T} \ge \tau \ge 0$.

PROPOSITION 4.2.   *The labeled sample partitions $\mathcal{P}^{\mathcal{C}_n}$ and $\mathcal{P}^{\mathcal{D}_n}$ generated by a marking with rate $\rho \, d\tau$ along $\mathcal{C}_n$ and $\mathcal{D}_n$, respectively, coincide with probability $1 - \mathcal{O}(\frac{\log \alpha}{\alpha})$.*

We need a lemma for the proof of this proposition. Denote by $T_\varepsilon^Z$ the time when $Z$ hits level $\varepsilon$ for the first time and denote by $T_c^{\mathcal{C}_n}$ the time when the number of lines in $\mathcal{C}_n$ increases from $n - 1$ to $n$.

LEMMA 4.3.   *Assume $\alpha$ is large enough such that $(\log \alpha)^2 \ge 2$. Let $\varepsilon(\alpha) = \frac{(\log \alpha)^2}{\alpha}$. Then*

$$(4.6) \qquad \mathbf{P}_0[T_{\varepsilon(\alpha)}^Z < T_c^{\mathcal{C}_n}] = \mathcal{O}\left(\frac{1}{(\log \alpha)^2}\right).$$

PROOF. Our proof rests on a Green's function calculation analogous to those in the proof of Proposition 3.4 and so, since we already have an expression for the Green's function for the process $X$, it is convenient to "undo" our time change. If we denote by $T_\varepsilon^X$ the time when $X$ hits level $\varepsilon$ for the first time and denote by $T_c^{\mathcal{C}_n}$ (resp. $T_c^{\mathcal{T}_n}$) the time of the first coalescence in $\mathcal{C}_n$ (resp. $\mathcal{T}_n$), then

$$\mathbf{P}_0[T_{\varepsilon(\alpha)}^Z < T_c^{\mathcal{C}_n}] = \mathbf{P}_0[T_{\varepsilon(\alpha)}^X < T_c^{\mathcal{T}_n}].$$

It is enough to consider the coalescence time of a sample of size 2, because as the probability that any pair in the sample coalesces is bounded by the sum over all pairs of lineages.

With (4.2),

$$\mathbf{P}_0[T_{\varepsilon(\alpha)}^X < T_c^{\mathcal{C}_2}] = 1 - \mathbf{E}_{\varepsilon(\alpha)}\left[\exp\left(-\int_0^T \frac{2}{X_s}\right)\right]$$

$$\leq \mathbf{E}_{\varepsilon(\alpha)}\left[\int_0^T \frac{2}{X_s}\,ds\right]$$

$$= \int_0^1 G(\varepsilon(\alpha), \xi)\frac{2}{\xi}\,d\xi.$$

Define $g(\alpha) := (\log \alpha)^2$. We split the last integral into three parts. We have, for constants $C$ which change from line to line,

$$\int_0^{\varepsilon(\alpha)} G(\varepsilon(\alpha), \xi)\frac{2}{\xi}\,d\xi = 2\int_0^{\varepsilon(\alpha)} \frac{(e^{-g(\alpha)} - e^{-\alpha})(e^{\alpha\xi} - 1)(1 - e^{-\alpha\xi})}{\alpha\xi^2(1-\xi)(1 - e^{-g(\alpha)})(1 - e^{-\alpha})}\,d\xi$$

$$\leq \frac{2e^{-g(\alpha)}}{(1 - e^{-g(\alpha)})(1 - \varepsilon(\alpha))}\int_0^{g(\alpha)} \frac{(e^\xi - 1)(1 - e^{-\xi})}{\xi^2}\,d\xi$$

$$\leq C\left(e^{-g(\alpha)} + e^{-g(\alpha)/2}\int_1^{g(\alpha)/2} d\xi + \int_{g(\alpha)/2}^{g(\alpha)} \frac{1}{\xi^2}\,d\xi\right)$$

$$\leq \frac{C}{g(\alpha)},$$

$$\int_{\varepsilon(\alpha)}^{1/2} G(\varepsilon(\alpha), \xi)\frac{2}{\xi}\,d\xi = 2\int_{\varepsilon(\alpha)}^{1/2} \frac{(1 - e^{-\alpha\xi})(1 - e^{-\alpha(1-\xi)})}{\alpha\xi^2(1-\xi)(1 - e^{-\alpha})}\,d\xi$$

$$\leq 4\int_{g(\alpha)}^{\alpha/2} \frac{1}{\xi^2}\,d\xi \leq \frac{4}{g(\alpha)},$$

$$\int_{1/2}^1 G(\varepsilon(\alpha), \xi)\frac{2}{\xi}\,d\xi \leq 4\int_{1/2}^1 G(0, \xi)\,d\xi \leq \frac{C\log \alpha}{\alpha}$$

as $g(\alpha) \geq 2$ and so $\varepsilon(\alpha) \leq \frac{1}{2}$, where we have used (3.1) in the third term. □

REMARK 4.4.   It is immediate from Lemma 4.3 that still writing $\varepsilon(\alpha) = \frac{(\log \alpha)^2}{\alpha}$ we will have

$$(4.7) \qquad \mathbf{P}_0[T^Z_{\varepsilon(\alpha)} < T^{\mathcal{D}_n}_c] = \mathcal{O}\left(\frac{1}{(\log \alpha)^2}\right).$$

PROOF OF PROPOSITION 4.2.   Let $T^Z_{\varepsilon(\alpha)}$ be as in Lemma 4.3. Looking backward from the time $T^Z_{\varepsilon(\alpha)}$, assume we take the marked $n$-coalescent with pair coalescence rate $\frac{2}{Z}$ as an approximation for the marked $n$-coalescent with pair coalescence rate $\frac{2}{Z(1-Z)}$. Sources of error in constructing the labeled partition are the recombination events that occur at a time when the two coalescents have different numbers of extant lines. We will call such an event a *bad* recombination event.

First we couple the two coalescent trees. We write $T^k_Z$, $T^k_{Z(1-Z)}$ for the times at which the coalescents with rates $2/Z$ and $2/Z(1-Z)$ per pair, respectively, have a transition from $k$ to $k-1$ extant lineages. We shall call our coupling successful if we have

$$T^n_{Z(1-Z)} > T^n_Z > T^{n-1}_{Z(1-Z)} > T^{n-1}_Z > \cdots > T^k_{Z(1-Z)} > T^k_Z > \cdots > T^2_{Z(1-Z)} > T^2_Z.$$

Let $S_1, \ldots, S_n$ be independent exponentially distributed random variables with $S_k \sim \exp(\binom{k}{2})$. The idea is simply to use the same random variable $S_k$ to generate the $k$th coalescence for both processes. Thus, writing $V$ for $Z$ or $Z(1-Z)$, $T^n_V, \ldots, T^1_V$ are defined recursively by

$$\int_{T^n_V}^{T^Z_{\varepsilon(\alpha)}} \frac{2}{V_s}\, ds = S_n \quad \text{and} \quad \int_{T^k_V}^{T^{k+1}_V} \frac{2}{V_s}\, ds = S_n, \qquad k = 1, \ldots, n-1.$$

Notice that our first inequality, $T^n_{Z(1-Z)} > T^n_Z$, is automatically satisfied. At time $T^n_{Z(1-Z)}$ we have

$$\int_{T^n_{Z(1-Z)}}^{T^Z_{\varepsilon(\alpha)}} \frac{2}{Z_s}\, ds \geq (1 - \varepsilon(\alpha)) S_n.$$

Thus

$$\begin{aligned}
\mathbf{P}[T^{n-1}_{Z(1-Z)} > T^n_Z] &= \mathbf{P}\left[\int_{T^{n-1}_{Z(1-Z)}}^{T^Z_{\varepsilon(\alpha)}} \frac{2}{Z_s}\, ds < \int_{T^{n-1}_Z}^{T^Z_{\varepsilon(\alpha)}} \frac{2}{Z_s}\, ds\right] \\
&\leq \mathbf{P}[(1 - \varepsilon(\alpha))(S_{n-1} + S_n) < S_n] \\
&= \mathcal{O}(\varepsilon(\alpha)).
\end{aligned}$$

Suppose then that $T^{n-1}_{Z(1-Z)} < T^n_Z$. Automatically then $T^{n-1}_Z < T^{n-1}_{Z(1-Z)}$ and at time $T^{n-1}_{Z(1-Z)}$, there is at most a further $\varepsilon(\alpha)(S_{n-1} + S_n)$ to accumulate in the integral defining $T_Z^{n-1}$ as

$$\mathbf{P}[T_{Z(1-Z)}^{n-2} > T_Z^{n-1}] = \mathbf{P}\left[\int_{T_{Z(1-Z)}^{n-2}}^{T_{\varepsilon(\alpha)}^Z} \frac{2}{Z_s}\,ds < \int_{T_Z^{n-1}}^{T_{\varepsilon(\alpha)}^Z} \frac{2}{Z_s}\,ds\right]$$

$$\leq \mathbf{P}[(1-\varepsilon(\alpha))(S_{n-2}+S_{n-1}+S_n) < S_{n-1}+S_n]$$

$$= \mathcal{O}(\varepsilon(\alpha)).$$

Continuing in this way we see that the chance of failing to achieve a successful coupling is $\mathcal{O}(\varepsilon(\alpha))$.

Now consider the chances of a *bad* recombination event on the tree when we have a successful coupling. We have to consider the lengths of the intervals when there are different numbers of lineages extant in the two coalescents, but these are the times $T_{Z(1-Z)}^k - T_Z^k$. We know that these are time intervals during which the $2/Z$ integral must accumulate on the order of $\varepsilon(\alpha)$. Since $Z < \varepsilon(\alpha)$, the time taken for this is $\mathcal{O}(\varepsilon(\alpha)^2)$. A *bad* recombination event then has probability $\mathcal{O}(\frac{\alpha}{\log \alpha}\varepsilon(\alpha)^2)$.

The errors that we have so far are the following:

- Failure to couple: an error of $\mathcal{O}(\varepsilon(\alpha))$.
- Coupling but bad recombinations: an error of $\mathcal{O}(\frac{\alpha}{\log \alpha}\varepsilon(\alpha)^2)$.

Since, by Lemma 4.3, the additional error coming from a coalescence of $\mathcal{C}_n$ or $\mathcal{D}_n$ between $T_{\varepsilon(\alpha)}^Z$ and $\tilde{T}$ is $\mathcal{O}((1/\log \alpha)^2)$, the proof of Proposition 4.2 is complete.  □

LEMMA 4.5.   *Let* $Z = (Z_\tau)_{0 \leq \tau < \infty}$ *be an* $\alpha$-*supercritical Feller process governed by*

$$dZ = \sqrt{2Z}\,dW + \alpha Z\,d\tau,$$

*started in* 0 *and conditioned on nonextinction, and let* $\mathcal{Y} = \mathcal{Y}^Z$ *be the tree of individuals with infinite lines of descent. Then the following statements are true:*

(a) *When averaged over* $Z$, $\mathcal{Y}$ *is a Yule tree with birth rate* $\alpha$.

(b) *The number of lines in* $\mathcal{Y}^Z$ *extant at time* $\tilde{T}$ *(the time when* $Z$ *first hits the level 1) has a Poisson distribution with mean* $\alpha$.

(c) *Given* $Z$, *the pair coalescence rate of* $\mathcal{Y}^Z$ *viewed backward from time* $\infty$ *is* $2/Z_\tau$.

PROOF.   Statement (a) is Theorem 3.2 of [19]. Statement (b) follows from [7] and the strong Markov property. Statement (c) derives from the fact that the pair coalescence rate of ancestral lineages in a continuum branching process (be it supercritical or not), conditioned on the total mass path $Z$, is a variant of Perkins' disintegration result, compare the discussion following Theorem 1.1 in [3].  □

REMARK 4.6. Part (c) of Lemma 4.5 also reveals that Kingman's coalescent with pair coalescence rate $\sigma^2/P$ describes the genealogy behind (1.1). Indeed, represent $P$ as

$$P_t = \frac{Z_\tau^{(1)}}{Z_\tau^{(1)} + Z_\tau^{(2)}},$$

where $Z^{(1)}$ is an $\alpha$-subcritical Feller process, $Z^{(2)}$ is a critical Feller process and the time change from $\tau$ to $t$ is given by

$$dt = \frac{d\tau}{Z_\tau^{(1)} + Z_\tau^{(2)}}.$$

Lemma 4.5(c) says that, conditional on $Z^{(1)}$, the pair coalescence rate in the genealogy of $Z^{(1)}$ is $\sigma^2 \, d\tau / Z_\tau^{(1)}$, which equals

$$\sigma^2 \, dt \frac{Z_\tau^{(1)} + Z_\tau^{(2)}}{Z_\tau^{(1)}} = \sigma^2 \, dt \frac{1}{P_t}.$$

PROPOSITION 4.7. *The variation distance between the distributions of the labeled partitions $\mathcal{P}^{\mathcal{Y}_n}$ (introduced in Definition 3.3) and $\mathcal{P}^{\mathcal{D}_n}$ is $\mathcal{O}(1/(\log \alpha)^2)$.*

PROOF. Our proof proceeds in two steps. First we account for the error that we make in assuming that there is no coalescence in the sample from the leaves of the Yule tree after time $\tilde{T}$. Second we account for the error in considering the process of marks on our Yule tree not up until time $\tilde{T}$ when there are a Poisson($\alpha$) number of extant individuals, but until the time when there are exactly $\lfloor \alpha \rfloor$ extant individuals.

(i) Given $Z$, take a sample of size $n$ from the leaves of $\mathcal{Y}^Z$. Write $\mathcal{Y}_n^Z$ for the ancestral tree of this sample and call $\mathcal{Y}_{n,\tilde{T}}^Z$ the cutoff of $\mathcal{Y}_n^Z$ between times $\tau = 0$ and $\tau = \tilde{T}$. Assume there are $N$ lines at time $\tilde{T}$. Since their lines of ascent agglomerate like in a Pólya urn, the proportions of their offspring in the leaves of $\mathcal{Y}^Z$ are uniformly distributed on the simplex $\{(p_1,\ldots,p_N)|p_h \geq 0, p_1 + \cdots + p_N = 1\}$. (See [10] and [13] for more background on Pólya urn schemes.) Writing $D_h$, $h = 1,\ldots,N$, for the number of all descendants of individual number $h$ which belong to the sample, one therefore obtains that $(D_1,\ldots,D_N)$ is uniformly distributed on

$$B_{N,n} := \{(d_1,\ldots,d_N)|d_1,\ldots,d_N \in \mathbb{N}_0, d_1 + \cdots + d_N = n\},$$

the set of occupation numbers of $N$ boxes with $n$ balls. (This distribution is also called the Bose–Einstein distribution with parameters $N$ and $n$.) Under this distribution the probability for multiple hits is

$$(4.8) \qquad 1 - \frac{\binom{N}{n}}{\binom{N+n-1}{n}} = \mathcal{O}\left(\frac{1}{N}\right) \qquad \text{as } N \to \infty.$$

Denote by $E_n$ the event that there is no coalescence of the ancestral lineages of the sample between times $\infty$ and $\tilde{T}$. Because of Lemma 4.5(b), the probability of $E_n$ arises by averaging (4.8) over a Poisson($\alpha$)-distributed $N$. Consequently the probability of $E_n$ when averaged over $Z$ is $\mathcal{O}(\frac{1}{\alpha})$.

(ii) Next we estimate the difference which it makes for the labeled partitions if we mark the branches of $\mathcal{Y}_n^Z$ at rate $\gamma\alpha/\log\alpha$ (1) between the real times $\tau = 0$ and $\tau = \tilde{T}$ or (2) between the "Yule times" $i = 1$ and $i = \lfloor\alpha\rfloor$. Since by Lemma 4.5(b) the number of lines extant in $\mathcal{Y}^Z$ at time $\tau = \tilde{T}$ is Poisson($\alpha$), to complete the proof of the proposition it suffices to show that the probability that $n$ chosen lines are hit by some mark between Yule times $\lfloor\alpha\rfloor$ and $J$ is $\mathcal{O}(1/(\log\alpha)^2)$, where $J$ has a Poisson($\alpha$) distribution.

Now on the one hand, by the Chebyshev inequality,

$$(4.9) \qquad \mathbf{P}[|J - \alpha| > \alpha^{3/4}] \leq \frac{\alpha}{(\alpha^{3/4})^2} = \alpha^{-1/2}.$$

On the other hand, from (3.6) we see that the probability that $n$ lines are hit by a mark between Yule times $i = \lfloor\alpha - \alpha^{3/4}\rfloor$ and $i = \lceil\alpha + \alpha^{3/4}\rceil$ is bounded by

$$1 - \exp\left(-\frac{n\gamma}{\log\alpha}\sum_{i=\lfloor\alpha-\alpha^{3/4}\rfloor}^{\lceil\alpha+\alpha^{3/4}\rceil}\frac{1}{i}\right) + \mathcal{O}\left(\frac{1}{(\log\alpha)^2}\right) \leq \frac{n\gamma}{\log\alpha}\frac{C\alpha^{3/4}}{\alpha} + \mathcal{O}\left(\frac{1}{(\log\alpha)^2}\right)$$

for a suitable $C > 0$.

Combining this with (4.9) and step (i), the assertion of Proposition 4.7 follows.  $\square$

Because of Corollary 3.5 and Propositions 4.2 and 4.7, and since the distributions of $\mathcal{P}^{\mathcal{T}_n}$ and $\mathcal{P}^{\mathcal{C}_n}$ coincide, Proposition 3.6 is now immediate. $\square$

### 4.4. *Within the Yule world*: *Proofs of Propositions* 3.8–3.11 *and* 2.8.

PROOF OF PROPOSITION 3.8.   All the results we are going to prove in this section deal with a sample of size $n$ taken from the leaves of an infinite Yule tree. As a key result we first obtain the "split times" in the sample genealogy as time evolves *forward* from $t = 0$. Recall from Section 3 that $I = I(t)$ is the number of lines of $\mathcal{Y}$ extant at time $t$ and that $K_i$ is the number of lines extant in $\mathcal{Y}_n$ while $I = i$.

LEMMA 4.8.   $K_i$, $i = 1, 2, \ldots$, *starts in* $K_1 = 1$ *and is a time-inhomogeneous Markov chain with transition probabilities*

$$(4.10) \quad \mathbf{P}[K_{i+1} = k + 1|K_i = k] = \frac{n - k}{n + i}, \qquad \mathbf{P}[K_{i+1} = k|K_i = k] = \frac{k + i}{n + i}.$$

*Its backward transition probabilities are*

$$(4.11) \qquad \mathbf{P}[K_i = k | K_{i+1} = k+1] = \frac{k(k+1)}{i(i+1)}.$$

*The one-time probabilities and the more-step forward and backward transition probabilities of $K$ are given by*

$$(4.12) \qquad \mathbf{P}[K_i = k] = \frac{\binom{n-1}{n-k}\binom{i}{k}}{\binom{n+i-1}{n}} \qquad (1 \le k \le \min(i,n)),$$

$$(4.13) \quad \mathbf{P}[K_j = l | K_i = k] = \frac{\binom{n-k}{n-l}\binom{j+k-1}{i+l-1}}{\binom{n+j-1}{n+i-1}},$$

$$(4.14) \quad \mathbf{P}[K_i = k | K_j = l] = \frac{\binom{j+k-1}{i+l-1}\binom{i}{k}\binom{l-1}{k-1}}{\binom{j-1}{i-1}\binom{j}{l}} \qquad (1 \le i, l \le j; k \le \min(i,l)).$$

PROOF. (i) We begin by deriving the one-step transition probabilities (4.10). At each Yule time $i$, attach to each line of $\mathcal{Y}$ the label $1+d$, where $d$ is the number of the line's descendants at infinity which belong to the sample. Call a line of $\mathcal{Y}$ *fertile* if it belongs to $\mathcal{Y}_n$, that is, if its attached label is larger than 1.

Passing from $i = 1$ to $i = 2$, this induces a split of the sample into subgroups of sizes $D_1$ and $D_2 = n - D_1$, where $D_1$ is uniform on $\{0, 1, \ldots, n\}$ (see the argument in the proof of Proposition 4.7). Given $(D_1, D_2)$, the proportion of the population at infinity in the tree $\mathcal{Y}$ that is descended from the line labeled $1 + D_1$ is Beta$(1 + D_1, 1 + n - D_1)$-distributed, the posterior of a Beta$(1,1)$ with $D_1$ successes in $n$ trials. Hence the birth in $\mathcal{Y}$ at Yule time $i = 2$ is to the line labeled $1 + D_1$ with probability $\frac{1+D_1}{n+2}$. If this is the case, $D_1$ is split uniformly into two subgroups, where a uniform split of 0 is understood as $(0,0)$.

At the $i$th stage of our construction there will be $i$ lines and an associated partition of $n + i$ into $i$ subsets with sizes denoted by

$$1 + D_1^i, 1 + D_2^i, \ldots, 1 + D_i^i,$$

where $D_1^i, \ldots, D_i^i$ are nonnegative integers that sum to $n$. The $(i+1)$st split is of the $j$th subset with probability $\frac{1+D_j}{i+n}$. At Yule time $i$, the probability of a "true split," that is, a split leading to two fertile successors of the fertile line labeled $1 + D_j$, is

$$\frac{1 + D_j}{i + n} \frac{D_j - 1}{D_j + 1} = \frac{D_j - 1}{i + n}.$$

Hence, given that the number $K_i$ of fertile lines at Yule time $i$ equals $k$, the probability of an increase of the number of fertile lines by 1 is

$$\frac{n-k}{i+n}.$$

Thus, the process $(K_i)_{i=1,2,\ldots}$ is a pure birth process in discrete time, starting in $K_1 = 1$, with time inhomogeneous transition probability given by (4.10).

Formula (4.11) follows from (4.10) and (4.12), which will be derived in the next step.

(ii) Next we derive the one-time probabilities (4.12). For this, take a sample of size $n$ from the leaves of the Yule tree $\mathcal{Y}$ and look at Yule time $i$, which is the period when there are $i$ individuals extant in $\mathcal{Y}$. Number these individuals by $h = 1, \ldots, i$ and let $D_h$, $h = 1, \ldots, i$, be the number of all descendants of individual number $h$ which belong to the sample. Then, by the argument given in the proof of Proposition 4.7, $(D_1, \ldots, D_i)$ is uniformly distributed on

$$B_{i,n} := \{(d_1, \ldots, d_i) | d_1, \ldots, d_i \in \mathbb{N}_0, d_1 + \cdots + d_i = n\},$$

the set of occupation numbers of $i$ boxes with $n$ balls. This distribution is also called the Bose–Einstein distribution with parameters $i$ and $n$.

The event "there are $k$ fertile lines at time $i$" thus has the same distribution as the event "$k$ of the $i$ boxes are occupied and the remaining $i-k$ are empty" under the Bose–Einstein distribution with parameters $i$ and $n$.

Let

$$B_{k,n}^+ := \{(d_1, \ldots, d_k) | d_1, \ldots, d_k \in \mathbb{N}, d_1 + \cdots + d_k = n\}.$$

It is easy to see that $\#B_{i,n} = \binom{n+i-1}{n}$ and $\#B_{k,n}^+ = \#B_{k,n-k} = \binom{n-1}{n-k}$. Hence, (4.12) arises as the probability under the uniform distribution on $B_{i,n}$ that $k$ of the $i$ boxes are occupied and the remaining $i-k$ are empty.

(iii) From (i) we see that our Markov chain can be represented in the following way: At Yule time $i$ there are $i$ *real* individuals in the Yule tree. Additionally, there are $n$ *virtual* individuals, corresponding to the sample of size $n$, each of which is attached to one of the real individuals. In this way, the $n$ virtual individuals are split in $k$ blocks; imagine that each block has one of its virtual individuals as its own *block leader*. Then evolve the Markov chain by choosing at each time one of the $n+i$ (real or virtual) individuals at random. When a real individual is chosen, it splits into two real individuals (leaving the number of fertile lines constant). When a virtual individual is chosen, then it is either a block leader or not. If it is a block leader, we find a split of a fertile line into one fertile and one infertile line (again leaving the number of blocks constant). When a non-block leader is picked, this gives rise to a new block and the chosen individual becomes a block leader. This then increases the number of blocks by 1.

We now proceed to prove (4.13). We start the Yule tree when it has $i$ lines, assuming there are currently $k$ blocks of virtual individuals. Thus there are currently $k$ block leaders and $n - k$ virtual individuals that are eligible to become new block leaders by time $j$. To distribute the additional $j - i$ individuals that enter the Yule tree between Yule times $i$ and $j$, note that these can choose among $i + n$ individuals as potential ancestors at time $i$. Thus the number of ways to distribute the $j - i$ additional individuals on $i + n$ ancestors is

$$\#B_{j-i,n+i} = \binom{n+j-1}{j-i} = \binom{n+j-1}{n+i-1}$$

and all of them have the same probability since the additional individuals arrive as in a Pólya scheme. To obtain $l$ blocks at Yule time $j$, having $k$ block leaders at Yule time $i$, we must calculate the probability of having added $l - k$ blocks up to Yule time $j$. At Yule time $i$, there are already $k$ block leaders and the number of ways to choose $l - k$ additional ones from the remaining $n - k$ potential new block leaders is

$$\binom{n-k}{l-k}.$$

To realize these new block leaders, we already have to use $l - k$ individuals. The remaining $j - i - (l - k)$ individuals must be distributed among the $i + n$ individuals present at Yule time $i$. However, to obtain $l$ blocks at Yule time $j$, these individuals must avoid the $n - l$ nonblock leaders (because this would result in new blocks). The number of ways to distribute $j - i - (l - k)$ balls on $i + n - (n - l)$ boxes is

$$\#B_{i+l,j-i-(l-k)} = \binom{j+k-1}{j-i-(l-k)} = \binom{j+k-1}{i+l-1}.$$

Altogether we arrive at (4.13).

(iv) The more-step backward transition probabilities (4.14) follow from (4.12), (4.13) and Bayes' formula.  □

REMARK 4.9. (i) Here is a self-contained derivation of the more-step backward transition probabilities (4.14). For $i, l \leq j$, $k \leq \min(i, l)$, consider the classical Pólya model with $i$ ancestors and $j - i$ newcomers, leading to $j$ people after the successive arrivals of the newcomers. Each of these $j - i$ newcomers joins the family of one of the $i$ ancestors by randomly choosing one of the extant individuals. The joint distribution of the numbers of newcomers in each family is then Bose–Einstein with parameters $i$ and $j - i$. There are now $j$ people in our population. Sample $l$ at random (without replacement). Then the conditional probability $\mathbf{P}[K_i = k | K_j = l]$ equals the probability that these $l$ people belong to exactly $k$ families. We can now decompose with respect to the number $U$ of individuals in the $l$-sample which are among the original $i$ ancestors. This number is hypergeometric, choosing $l$ out of $j = i + (j - i)$ [see (2.6)].

Given $U = u$, the probability that the sample forms exactly $k$ distinct families is the probability that exactly $k - u$ of the available $i - u$ ancestors have descendants among the $l - u$ newcomers in the sample. There are $\binom{i-u}{k-u}$ possible choices of the additional ancestors and, conditional on that choice, there are $\binom{l-u-(k-u)+k-1}{l-u-(k-u)}$ ways to distribute the remaining newcomers. Hence the conditional probability that the sample forms exactly $k$ distinct families is

$$\frac{\binom{i-u}{k-u}\binom{l-u-(k-u)+k-1}{l-u-(k-u)}}{\binom{l-u+i-1}{l-u}} = \frac{\binom{i-u}{i-k}\binom{l-1}{l-k}}{\binom{l-u+i-1}{i-1}}.$$

We thus conclude that

$$(4.15) \qquad \mathbf{P}[K_i = k | K_j = l] = \sum_{u=0}^{k} \frac{\binom{i}{u}\binom{j-i}{l-u}}{\binom{j}{l}} \frac{\binom{i-u}{i-k}\binom{l-1}{l-k}}{\binom{l-u+i-1}{i-1}}.$$

The fact that the right-hand sides of (4.14) and (4.15) are equal can be checked by an elementary but tedious calculation (or by your favorite computer algebra package).

(ii) Since (4.12) follows from (4.14) by putting $l = n$ and letting $j$ tend to infinity, (4.13) follows from (4.14) and (4.12) by the Bayes formula, and (4.10) and (4.11) specialize from (4.13) and (4.14).

With Lemma 4.8 it is now easy to complete the proof of Proposition 3.8:

$$\mathbf{P}[F \leq i] = \mathbf{P}[K_i = n]$$

$$= \frac{\binom{i}{n}}{\binom{n+i-1}{n}} = \frac{i!(i-1)!}{(i-n)!(n+i-1)!}$$

$$= \frac{(i-1)\cdots(i-n+1)}{(i+n-1)\cdots\cdots(i+1)}. \qquad \qquad \square$$

**PROOF OF PROPOSITION 3.9.** As a preparation, we prove the following lemma:

**LEMMA 4.10.** *For the first Yule time $F$ when there are $n$ fertile lines we have, for $k < n$,*

$$\mathbf{P}[F = f | K_i = k] = \frac{(f-i-1)\cdots\cdots(f-i-(n-k)+1)}{(f+n-1)\cdots\cdots(f+k-1)}$$

$$(4.16) \qquad\qquad\qquad \times (n-k)(n+i-1), \qquad k < n-1,$$

$$\mathbf{P}[F = f | K_i = n-1] = \frac{1}{(f+n-1)(f+n-2)}(n+i-1)$$

*and*

(4.17)
$$\mathbf{P}[F = f | K_i = k] \le \frac{Ci}{f^2},$$

*where $C$ depends only on $n$.*

PROOF.  Equation (4.16) follows by

$$\mathbf{P}[F = f | K_i = k]$$
$$= \mathbf{P}[F > f - 1 | K_i = k] - \mathbf{P}[F > f | K_i = k]$$
$$= \mathbf{P}[K_f = n | K_i = k] - \mathbf{P}[K_{f-1} = n | K_i = k]$$
$$= \frac{\binom{f+k-1}{i+n-1}}{\binom{n+f-1}{n+i-1}} - \frac{\binom{f+k-2}{i+n-1}}{\binom{n+f-2}{n+i-1}}$$
$$= \frac{(f+k-1)!(f-i)!}{(f-i+k-n)!(n+f-1)!} - \frac{(f+k-2)!(f-i-1)!}{(f-i+k-n-1)!(n+f-2)!}$$
$$= \{(f+k-2)!(f-i-1)!$$
$$\quad \times ((f+k-1)(f-i) - (f-i+k-n)(n+f-1))\}$$
$$\quad \times \{(f-i+k-n)!(n+f-1)!\}^{-1}$$
$$= \frac{(f+k-2)!(f-i-1)!(n-k)(n+i-1)}{(f-i+k-n)!(n+f-1)!},$$

where we have used

$$(f+k-1)(f-i) - (f-i+k-n)(n+f-1)$$
$$= (f-i)(k-n) - (k-n)(n+f-1)$$
$$= (k-n)(f-i-n-f+1) = (n-k)(n+i-1).$$

If $k < n - 1$, the terms $(f-i-1)!$ and $(f-i+k-n)!$ cancel partially, leading to

$$\mathbf{P}[F = f | K_i = k] = \frac{(f-i-1) \cdot \cdots \cdot (f-i-(n-k)+1)}{(f+n-1) \cdot \cdots \cdot (f+k-1)}(n-k)(n+i-1).$$

In the case $k = n - 1$, these two terms cancel completely, which proves (4.16).

To see (4.17), note that the fraction in (4.16) is bounded by $1/f^2$ and $(n-k)(n+i-1) \le n^2 i$.  □

Recall the definitions of $M_i$, $M = \sum_{i=1}^{F} M_i$ and $S^{\mathcal{Y}}$ from Section 3. We first turn to the proof of (3.5).

First observe that since $\{M \geq 2\}$ requires either that the tree is hit by marks during two distinct Yule times in the early phase or at least twice during a single Yule time (during the early phase), we may estimate

$$\mathbf{P}[M \geq 2] \leq \sum_{i=n}^{\lfloor \alpha \rfloor} \sum_{j=i+1}^{\lfloor \alpha \rfloor} \mathbf{P}[M_i \geq 1, M_j \geq 1, K_j < n] + \sum_{i=n}^{\lfloor \alpha \rfloor} \mathbf{P}[M_i \geq 2].$$

For the first term we use (4.17) to see that there is some constant $C$ that depends only on $n$ such that, for $k < n$,

$$\mathbf{P}[F > j | K_i = k] \leq C \frac{i}{j}.$$

With this we approximate (recall Remark 3.7) for constants $C$ which depend only on $n$ and can change from occurrence to occurrence:

$$\sum_{i=1}^{\lfloor \alpha \rfloor} \sum_{j=i+1}^{\lfloor \alpha \rfloor} \mathbf{P}[M_i \geq 1, M_j \geq 1, K_j < n]$$

$$\leq \sum_{i=1}^{\lfloor \alpha \rfloor} \sum_{j=i+1}^{\lfloor \alpha \rfloor} \sum_{k_1, k_2 = 1}^{n-1} \mathbf{P}[M_j = 1 | K_j = k_2] \cdot \mathbf{P}[M_i = 1 | K_i = k_1]$$

$$\times \mathbf{P}[K_i = k_1, K_j = k_2]$$

$$\leq \sum_{i=1}^{\lfloor \alpha \rfloor} \sum_{j=i+1}^{\lfloor \alpha \rfloor} \sum_{k_1, k_2 = 1}^{n-1} \frac{\gamma^2 n^2}{(\log \alpha)^2} \frac{1}{ij} \mathbf{P}[K_i = k_1, K_j = k_2]$$

$$= \frac{\gamma^2 n^2}{(\log \alpha)^2} \sum_{i=1}^{\lfloor \alpha \rfloor} \sum_{j=i+1}^{\lfloor \alpha \rfloor} \sum_{k_1 = 1}^{n-1} \frac{1}{ij} \mathbf{P}[K_j \leq n - 1 | K_i = k_1] \cdot \mathbf{P}[K_i = k_1]$$

$$\leq \frac{C}{(\log \alpha)^2} \sum_{i=1}^{\lfloor \alpha \rfloor} \sum_{j=i+1}^{\lfloor \alpha \rfloor} \frac{1}{ij} \sum_{k_1 = 1}^{n-1} \mathbf{P}[F > j | K_i = k_1] \cdot \mathbf{P}[K_i = k_1]$$

$$\leq \frac{C}{(\log \alpha)^2} \sum_{i=1}^{\lfloor \alpha \rfloor} \sum_{j=i+1}^{\lfloor \alpha \rfloor} \frac{1}{ij} \frac{i}{j} \mathbf{P}[F > i]$$

$$\leq \frac{C}{(\log \alpha)^2} \sum_{i=1}^{\lfloor \alpha \rfloor} \sum_{j=i+1}^{\lfloor \alpha \rfloor} \frac{1}{ij^2}$$

$$\leq \frac{C}{(\log \alpha)^2} \sum_{i=1}^{\lfloor \alpha \rfloor} \frac{1}{i^2} \leq \frac{C}{(\log \alpha)^2}.$$

For the second term, again by Remark 3.7,

$$\sum_{i=1}^{\lfloor \alpha \rfloor} \mathbf{P}[M_i \geq 2] \leq \sum_{i=1}^{\lfloor \alpha \rfloor} \frac{C}{i^2 (\log \alpha)^2} \leq \frac{C}{(\log \alpha)^2},$$

which proves (3.5).

Consequently, we can now concentrate on the event $\{M = 1\}$. The proof of (3.4) consists of two steps. First we will show that

$$
\begin{aligned}
(4.18) \quad & \mathbf{P}[S^{\mathcal{Y}} = s, M = 1] \\
& = \mathcal{O}\left(\frac{1}{(\log \alpha)^2}\right) + \begin{cases} \dfrac{\gamma}{\log \alpha} \displaystyle\sum_{i=1}^{\lfloor \alpha \rfloor} \dfrac{\binom{n-s+i-2}{n-s}}{\binom{n+i-1}{n}}, & s \geq 2, \\[2ex] \dfrac{\gamma}{\log \alpha} \displaystyle\sum_{i=1}^{\lfloor \alpha \rfloor} \dfrac{\binom{n+i-3}{n-1} - \binom{i-1}{n-1}}{\binom{n+i-1}{n}}, & s = 1. \end{cases}
\end{aligned}
$$

Second we will approximate these probabilities to obtain (3.4).

Given that the sample genealogy is hit by exactly one mark, and given that this happens when there are $k$ fertile lines ($1 \leq k \leq n-1$), then the (conditional) probability that $S^{\mathcal{Y}} = s$ ($1 \leq s \leq n-k+1$) is (in the notation of the proof of Lemma 4.8)

$$
\begin{aligned}
(4.19) \quad \mathbf{P}[S^{\mathcal{Y}} = s | K_i = k, M = M_i = 1] &= \frac{\#B_{k-1, n-k-(s-1)}}{\#B_{k, n-k}} \\
&= \frac{\binom{n-k-(s-1)+(k-1)-1}{n-k-(s-1)}}{\binom{n-k+k-1}{n-k}} = \frac{\binom{n-s-1}{n-s-(k-1)}}{\binom{n-1}{n-k}}.
\end{aligned}
$$

[Indeed, the one line which is hit by a mark during period $i$ must spawn $s-1$ offspring and the other $k-1$ lines must spawn $n-k-(s-1)$ offspring.]

The probability that the sample genealogy is hit by a mark during period $i$ *and* there are $k$ fertile lines [$k \leq i \wedge (n-1)$] during period $i$ is, according to Remark 3.7 and (3.6),

$$
\begin{aligned}
(4.20) \quad \mathbf{P}[M_i = 1, K_i = k] &= \left(\frac{\gamma}{\log \alpha} \frac{k}{i} + \mathcal{O}\left(\frac{1}{i^2 (\log \alpha)^2}\right)\right) \mathbf{P}[K_i = k] \\
&= \frac{\gamma}{\log \alpha} \frac{k}{i} \frac{\binom{i}{k}\binom{n-1}{n-k}}{\binom{n+i-1}{n}} + \mathcal{O}\left(\frac{1}{i^2 (\log \alpha)^2}\right) \\
&= \frac{\gamma}{\log \alpha} \frac{\binom{i-1}{k-1}\binom{n-1}{n-k}}{\binom{n+i-1}{n}} + \mathcal{O}\left(\frac{1}{i^2 (\log \alpha)^2}\right).
\end{aligned}
$$

To now show (4.18), we need some approximations.

LEMMA 4.11.  *For constants $C_1$ and $C_2$ which only depend on $n$,*

$$(4.21) \qquad \sum_{i=1}^{\lfloor \alpha \rfloor} \sum_{k=1}^{i \wedge (n-1)} (\log i) \mathbf{P}[K_i = k, M_i = 1] \leq \frac{C_1}{\log \alpha},$$

$$(4.22) \qquad \mathbf{P}[M \geq 2 | M_i = 1, K_i = k] \leq \frac{C_2(1 + \log i)}{\log \alpha}.$$

PROOF.  We will use the integral

$$(4.23) \qquad \int \frac{\log x}{x^2} = -\frac{1 + \log x}{x},$$

which will give us finiteness of some constants.

For (4.21), we have, from (4.20),

$$\sum_{i=1}^{\lfloor \alpha \rfloor} \sum_{k=1}^{i \wedge (n-1)} (\log i) \mathbf{P}[K_i = k, M_i = 1]$$

$$= \mathcal{O}\left(\frac{1}{(\log \alpha)^2}\right) + \frac{\gamma}{\log \alpha} \sum_{i=1}^{\lfloor \alpha \rfloor} \log i \sum_{k=1}^{i \wedge (n-1)} \frac{\binom{i-1}{k-1}\binom{n-1}{n-k}}{\binom{n+i-1}{n}}$$

$$= \mathcal{O}\left(\frac{1}{(\log \alpha)^2}\right) + \frac{\gamma}{\log \alpha} \sum_{i=1}^{\lfloor \alpha \rfloor} \log i \frac{\binom{n+i-2}{n-1} - \binom{i-1}{n-1}}{\binom{n+i-1}{n}}.$$

Now, for some constants $C$ and $C'$ which are bounded in $i$ (and may change from appearance to appearance),

$$\binom{n+i-2}{n-1} - \binom{i-1}{n-1} = \frac{[(i+n-2)\cdots i] - [(i-1)\cdots(i-(n-1))]}{(n-1)!}$$

$$\leq \frac{i^{n-1} + Ci^{n-2} - i^{n-1} + C'i^{n-2}}{(n-1)!} \leq Ci^{n-2}$$

and, because for $i \geq 2$

$$\binom{n+i-1}{n} \geq Ci^n,$$

we see that

$$\sum_{i=1}^{\lfloor \alpha \rfloor} \sum_{k=1}^{i \wedge (n-1)} (\log i) \mathbf{P}[K_i = k, M = M_i = 1]$$

$$\leq \mathcal{O}\left(\frac{1}{(\log \alpha)^2}\right) + \frac{\gamma}{\log \alpha} \sum_{i=1}^{\lfloor \alpha \rfloor} \frac{Ci^{n-2} \log i}{Ci^n} \leq \frac{C'}{\log \alpha},$$

where we have used (4.23).

For (4.22), using (4.17) and (4.23), we write, for some constant $C$ depending only on $n$ which can change from instance to instance,

$$\mathbf{P}[M \geq 2 | M_i = 1, K_i = k]$$

$$(4.24) \qquad = \sum_{j=i}^{\lfloor \alpha \rfloor} \mathbf{P}[M \geq 2 | M_i = 1, K_i = k, F = j] \cdot \mathbf{P}[F = j | K_i = k]$$

$$\leq \sum_{j=i}^{\lfloor \alpha \rfloor} \sum_{l=1}^{j} \frac{n}{l} \frac{\gamma}{\log \alpha} \frac{Ci}{j^2} \leq \frac{C}{\log \alpha} \sum_{j=i}^{\lfloor \alpha \rfloor} \frac{i \log j}{j^2} \leq \frac{C(1 + \log i)}{\log \alpha}. \qquad \square$$

We now return to the proof of (4.18). The probability that the sample genealogy is hit by a unique early mark *and* that $S^{\mathcal{Y}} = s \leq n$ is

$$\mathbf{P}[S^{\mathcal{Y}} = s, M = 1]$$

$$(4.25) \qquad = \sum_{k=1}^{i \wedge (n-1)} \sum_{i=1}^{\lfloor \alpha \rfloor} \mathbf{P}[S^{\mathcal{Y}} = s, M = M_i = 1, K_i = k]$$

$$= \sum_{k=1}^{i \wedge (n-1)} \sum_{i=1}^{\lfloor \alpha \rfloor} \mathbf{P}[S^{\mathcal{Y}} = s | M = M_i = 1, K_i = k]$$

$$(4.26) \qquad \qquad \times \mathbf{P}[M_i = 1, K_i = k](1 - \mathbf{P}[M \geq 2 | M_i = 1, K_i = k]).$$

For this sum, the event $\{M \geq 2\}$ does not play a role, because by (4.21) and (4.22),

$$\sum_{k=1}^{i \wedge (n-1)} \sum_{i=1}^{\lfloor \alpha \rfloor} \mathbf{P}[M_i = 1, K_i = k] \cdot \mathbf{P}[M \geq 2 | M_i = 1, K_i = k]$$

$$(4.27) \qquad \leq \sum_{i=1}^{\lfloor \alpha \rfloor} \frac{C \log i}{\log \alpha} \mathbf{P}[M_i = 1, K_i \leq n - 1] \leq \frac{C}{(\log \alpha)^2}.$$

So, combining (4.19) and (4.20), we have

$$\mathbf{P}[S^{\mathcal{Y}} = s, M = 1]$$

$$= \mathcal{O}\left(\frac{1}{(\log \alpha)^2}\right) + \sum_{k=1}^{i \wedge (n-1)} \sum_{i=1}^{\lfloor \alpha \rfloor} \mathbf{P}[S^{\mathcal{Y}} = s | M = M_i = 1, K_i = k]$$

$$\times \mathbf{P}[M_i = 1, K_i = k]$$

$$(4.28) \qquad = \mathcal{O}\left(\frac{1}{(\log \alpha)^2}\right) + \frac{\gamma}{\log \alpha} \sum_{i=1}^{\lfloor \alpha \rfloor} \sum_{k=1}^{i \wedge (n-1)} \frac{\binom{n-s-1}{n-s-(k-1)}\binom{i-1}{k-1}}{\binom{n+i-1}{n}}$$

$$= \mathcal{O}\Big(\frac{1}{(\log\alpha)^2}\Big) + \frac{\gamma}{\log\alpha}\sum_{i=1}^{\lfloor\alpha\rfloor}\frac{\binom{n-s+i-2}{n-s} - \binom{i-1}{n-1}\binom{n-s-1}{1-s}}{\binom{n+i-1}{n}}$$

$$= \mathcal{O}\Big(\frac{1}{(\log\alpha)^2}\Big) + \begin{cases} \dfrac{\gamma}{\log\alpha}\sum_{i=1}^{\lfloor\alpha\rfloor}\dfrac{\binom{n-s+i-2}{n-s}}{\binom{n+i-1}{n}}, & s\geq 2, \\[3mm] \dfrac{\gamma}{\log\alpha}\sum_{i=1}^{\lfloor\alpha\rfloor}\dfrac{\binom{n+i-3}{n-1} - \binom{i-1}{n-1}}{\binom{n+i-1}{n}}, & s=1, \end{cases}$$

which proves (4.18).

We now approximate these probabilities further to obtain (3.4). First observe that, since

$$\frac{1}{(n+i-1)\cdots i} = \frac{1}{n-1}\Big(\frac{1}{(n+i-2)\cdots i} - \frac{1}{(n+i-1)\cdots(i+1)}\Big)$$

for $s = n$ and $n \geq 2$, we may write

$$(4.29) \qquad \begin{aligned} \sum_{i=1}^{\lfloor\alpha\rfloor}\frac{1}{\binom{n+i-1}{n}} &= n!\sum_{i=1}^{\lfloor\alpha\rfloor}\frac{1}{(n+i-1)\cdots\cdot i} \\ &= \frac{n!}{(n-1)\cdot(n-1)!} + \mathcal{O}\Big(\frac{1}{\alpha}\Big) = \frac{n}{n-1} + \mathcal{O}\Big(\frac{1}{\alpha}\Big), \end{aligned}$$

which gives (3.4) in the case $s = n$.

Now define

$$(4.30) \qquad A(n,s,\alpha) := \sum_{i=1}^{\lfloor\alpha\rfloor}\frac{\binom{n-s+i-2}{n-s}}{\binom{n+i-1}{n}}.$$

For $s \leq n-1$, the summand vanishes for $i = 1$ and so

$$(4.31) \qquad \begin{aligned} A(n,s,\alpha) &= \sum_{i=1}^{\lfloor\alpha\rfloor-1}\frac{\binom{n-s+i-1}{n-s}}{\binom{n+i}{n}} = \sum_{i=1}^{\lfloor\alpha\rfloor-1}\frac{n!\,i!\,(n-s+i-1)!}{(n+i)!\,(n-s)!\,(i-1)!} \\ &= \frac{n!}{(n-s)!}\sum_{i=1}^{\lfloor\alpha\rfloor-1}\frac{i}{(i+n)\cdots(i+n-s)} \\ &= \frac{n!}{(n-s)!}\sum_{i=1}^{\lfloor\alpha\rfloor-1}\Big(\frac{1}{(i+n)\cdots(i+n-s+1)} \\ &\qquad\qquad - \frac{n-s}{(i+n)\cdots(i+n-s)}\Big). \end{aligned}$$

We treat the two sums separately and rewrite each as a telescoping sum as in the derivation of (4.29) to see that, for $2 \le s \le n-1$, this gives

$$
\begin{aligned}
A(n, s, \alpha) &= \frac{n!}{(n-s)!} \left( \frac{1}{s-1} \frac{1}{n \cdots (n-s+2)} - \frac{n-s}{s} \frac{1}{n \cdots (n-s+1)} \right) \\
&\quad + \mathcal{O}\left( \frac{1}{\alpha} \right) \\
&= \frac{n-s+1}{s-1} - \frac{n-s}{s} + \mathcal{O}\left( \frac{1}{\alpha} \right) \\
&= \frac{(n-s+1)s - (s-1)(n-s)}{(s-1)s} + \mathcal{O}\left( \frac{1}{\alpha} \right) \\
&= \frac{n}{(s-1)s} + \mathcal{O}\left( \frac{1}{\alpha} \right),
\end{aligned}
$$

so (3.4) is also proved for $2 \le s \le n-1$. For $s = 1$, the above gives

$$
\begin{aligned}
A(n, 1, \alpha) &= n \sum_{i=1}^{\lfloor \alpha \rfloor - 1} \frac{1}{i+n} - n(n-1) \sum_{i=1}^{\lfloor \alpha \rfloor - 1} \left( \frac{1}{i+n-1} - \frac{1}{i+n} \right) \\
&= n \sum_{i=n+1}^{\lfloor \alpha \rfloor} \frac{1}{i} - n(n-1) \frac{1}{n} + \mathcal{O}\left( \frac{1}{\alpha} \right) \\
&= 1 - n + n \sum_{i=n+1}^{\lfloor \alpha \rfloor} \frac{1}{i} + \mathcal{O}\left( \frac{1}{\alpha} \right).
\end{aligned}
\tag{4.32}
$$

For $s = 1$, we also have to deal with the second term in (4.18). We write

$$
\sum_{i=1}^{\lfloor \alpha \rfloor} \frac{\binom{i-1}{n-1}}{\binom{n+i-1}{n}} = \sum_{i=1}^{\lfloor \alpha \rfloor} \frac{(i-1)!(i-1)!n!}{(n-1)!(i-n)!(n+i-1)!} = n \sum_{i=1}^{\lfloor \alpha \rfloor} \frac{(i-1) \cdots (i-n+1)}{(i+n-1) \cdots i}.
$$

Define

$$
A_{m,n} := \sum_{i=1}^{\lfloor \alpha \rfloor} \frac{(i-1) \cdots (i-m+1)}{(i+n-1) \cdots i}, \qquad m > 1,
$$

$$
A_{1,n} := \sum_{i=1}^{\lfloor \alpha \rfloor} \frac{1}{(i+n-1) \cdots i}, \qquad A_n := A_{n,n}.
$$

Our goal then is to find an approximation of $A_n$. Observe that we have a recursive structure:

$$
A_{m,n} = \sum_{i=1}^{\lfloor \alpha \rfloor} \frac{(i-1) \cdots (i-m+2)}{(i+n-2) \cdots i}
$$

$$(4.33) \qquad -(m+n-2)\sum_{i=1}^{\lfloor \alpha \rfloor}\frac{(i-1)\cdots(i-m+2)}{(i+n-1)\cdots i}$$

$$= A_{m-1,n-1} - (m+n-2)A_{m-1,n}.$$

From this equation it also follows that

$$A_{m,n} = A_{1,n-m+1} + \sum_{k=0}^{m-2}(A_{m-k,n-k} - A_{m-k-1,n-k-1})$$

$$(4.34)$$

$$= A_{1,n-m+1} - \sum_{k=0}^{m-2}(m+n-2k-2)A_{m-1-k,n-k}.$$

First we show that for $1 \le m < n$,

$$(4.35) \qquad A_{m,n} = \frac{(m-1)!(m-1)!(n-m-1)!}{(n-1)!(n-1)!} + \mathcal{O}\Big(\frac{1}{\alpha}\Big).$$

We proceed by induction on $m$. For $m=1$, we have, up to an error of $\mathcal{O}(\frac{1}{\alpha})$ (here $n > 1$ is important),

$$A_{1,n} = \sum_{i=1}^{\lfloor \alpha \rfloor}\frac{1}{(i+n-1)\cdots i}$$

$$= \frac{1}{n-1}\sum_{i=1}^{\lfloor \alpha \rfloor}\frac{1}{(i+n-2)\cdots i} - \frac{1}{(i+n-1)\cdots(i+1)}$$

$$= \frac{1}{(n-1)(n-1)!} + \mathcal{O}\Big(\frac{1}{\alpha}\Big)$$

$$= \frac{(n-2)!}{(n-1)!(n-1)!} + \mathcal{O}\Big(\frac{1}{\alpha}\Big)$$

and, by (4.33),

$$A_{m+1,n} = A_{m,n-1} - (m+n-1)A_{m,n}$$

$$= \frac{(m-1)!(m-1)!(n-m-2)!}{(n-2)!(n-2)!}$$

$$\qquad -(m+n-1)\frac{(m-1)!(m-1)!(n-m-1)!}{(n-1)!(n-1)!}$$

$$= \frac{(m-1)!(m-1)!(n-m-2)!((n-1)^2 - (n-1+m)(n-1-m))}{(n-1)!(n-1)!}$$

$$= \frac{m!m!(n-(m+1)-1)!}{(n-1)!(n-1)!},$$

which proves (4.35).

From (4.35) and (4.34), we see, because

$$(4.36) \qquad A_1 = \sum_{i=1}^{\lfloor \alpha \rfloor} \frac{1}{i} + \mathcal{O}\left(\frac{1}{\alpha}\right),$$

that

$$
\begin{aligned}
(4.37) \quad A_n &= A_1 - 2\sum_{k=0}^{n-2}(n-k-1)\frac{(n-k-2)!(n-k-2)!}{(n-k-1)!(n-k-1)!} + \mathcal{O}\left(\frac{1}{\alpha}\right) \\
&= A_1 - 2\sum_{k=0}^{n-2}\frac{1}{n-k-1} + \mathcal{O}\left(\frac{1}{\alpha}\right) = A_1 - 2\sum_{i=1}^{n-1}\frac{1}{i} + \mathcal{O}\left(\frac{1}{\alpha}\right) \\
&= \frac{1}{n} - 1 + \sum_{i=n+1}^{\lfloor \alpha \rfloor}\frac{1}{i} - \sum_{i=2}^{n-1}\frac{1}{i} + \mathcal{O}\left(\frac{1}{\alpha}\right).
\end{aligned}
$$

Now (3.4) follows in the case $s=1$ from (4.18), (4.32) and (4.37) as

$$
\begin{aligned}
\sum_{i=1}^{\lfloor \alpha \rfloor} \frac{\binom{n+i-3}{n-1}-\binom{i-1}{n-1}}{\binom{n+i-1}{n}} &= 1 - n + n\sum_{i=n+1}^{\lfloor \alpha \rfloor}\frac{1}{i} \\
&\quad - n\left(\frac{1}{n} - 1 + \sum_{i=n+1}^{\lfloor \alpha \rfloor}\frac{1}{i} - \sum_{i=2}^{n-1}\frac{1}{i}\right) + \mathcal{O}\left(\frac{1}{\alpha}\right) \\
&= n\sum_{i=2}^{n-1}\frac{1}{i} + \mathcal{O}\left(\frac{1}{\alpha}\right). \qquad \square
\end{aligned}
$$

PROOF OF PROPOSITION 3.10.   Let $\mathcal{L}$ denote the random subset that consists of all those ancestral lineages of the sample which are hit by a late mark. For a fixed subset $A$ of the $n$ ancestral lineages of the sample, we conclude from Remark 3.7 and (3.6) that, for all $f \in \{1, \ldots, \lfloor \alpha \rfloor\}$,

$$
\begin{aligned}
(4.38) \quad \mathbf{P}[\mathcal{L} \cap A = \varnothing | F = f] &= \prod_{i=f}^{\lfloor \alpha \rfloor}\frac{i\alpha}{i\alpha + a\gamma\alpha/\log\alpha} \\
&= \exp\left(-\frac{a\gamma}{\log\alpha}\sum_{i=f}^{\lfloor \alpha \rfloor}\frac{1}{i}\right) + \mathcal{O}\left(\frac{1}{(\log\alpha)^2}\right) \\
&= (p_f)^a + \mathcal{O}\left(\frac{1}{(\log\alpha)^2}\right),
\end{aligned}
$$

where $a = \#A$ and the error term is uniform in $f$. Consequently, if we consider the random subset $\mathcal{M}$ of the sample which results from the successes of coin tossing with random success probability $p_F$, we observe that

$$(4.39) \quad \mathbf{P}[\mathcal{L} \cap A = \varnothing] = \mathbf{E}[(p_F)^a] + \mathcal{O}\left(\frac{1}{(\log \alpha)^2}\right)$$

$$= \mathbf{P}[\mathcal{M} \cap A = \varnothing] + \mathcal{O}\left(\frac{1}{(\log \alpha)^2}\right).$$

By inclusion–exclusion, (4.38) extends to the desired approximate equality of the distributions of $\mathcal{L}$ and $\mathcal{M}$.  $\square$

PROOF OF PROPOSITION 3.11. It remains to prove the approximate independence of the random variables $S^{\mathcal{Y}}$ and $L^{\mathcal{Y}}$. It is enough to show that, with the desired accuracy, $S^{\mathcal{Y}}$ is independent of the event of a late recombination of a randomly chosen line.

For convenience, we abuse notation, and write $S$ and $L$ instead of $S^{\mathcal{Y}}$ and $L^{\mathcal{Y}}$. The approximate independence of the distributions of $S$ and $L$ relies on two crucial observations. First, because $S > 0$ only has a probability of order $\mathcal{O}(\frac{1}{\log \alpha})$ (see Proposition 3.9), we can allow for a multiplicative error of order $\mathcal{O}(\frac{1}{\log \alpha})$ for the probability of $L = l$. The second observation is that the two probabilities $\mathbf{P}[L = l]$ and $\mathbf{P}[L = l | K_i = k]$ are $\frac{C \log i}{\log \alpha}$ apart, which is the content of (4.41).

LEMMA 4.12. *There are constants $C_1$ and $C_2$ that depend only on $n$ such that*

$$(4.40) \quad |\mathbf{P}[L = l | F = f] - \mathbf{P}[L = l | F = f']| \leq \frac{C_1(1 + \log f)(1 + \log f')}{\log \alpha},$$

$$(4.41) \quad |\mathbf{P}[L = l | K_i = k] - \mathbf{P}[L = l]| \leq \frac{C_2(1 + \log i)}{\log \alpha} \qquad (k < n).$$

PROOF. We start by proving (4.40). Given $F = f$, the number of late recombinants is approximately binomially distributed with parameters $n$ and $1 - p_f$ (see Proposition 3.10). Thus, for $f, f' \leq \alpha$,

$$|\mathbf{P}[L = l | F = f] - \mathbf{P}[L = l | F = f']|$$

$$= \binom{n}{l} |(1 - p_f)^l p_f^{n-l} - (1 - p_{f'})^l p_{f'}^{n-l}| + \mathcal{O}\left(\frac{1}{(\log \alpha)^2}\right)$$

$$\leq \sum_{k=0}^{l} \binom{n}{l} \binom{l}{k} |p_f^{n-l+k} - p_{f'}^{n-l+k}| + \mathcal{O}\left(\frac{1}{(\log \alpha)^2}\right),$$

where we have used that $(1-p)^l = \sum_{k=0}^l \binom{l}{k} p^k$. Now (4.40) follows, since for $0 \le m \le 2n$ and (w.l.o.g.) $f \le f'$,

$$
\begin{aligned}
|p_f^m - p_{f'}^m| &= \left| \exp\left( -\frac{m\gamma}{\log \alpha} \sum_{i=f}^{\lfloor \alpha \rfloor} \frac{1}{i} \right) - \exp\left( -\frac{m\gamma}{\log \alpha} \sum_{i=f'}^{\lfloor \alpha \rfloor} \frac{1}{i} \right) \right| \\
&\le \left| 1 - \exp\left( -\frac{m\gamma}{\log \alpha} \sum_{i=f}^{f'} \frac{1}{i} \right) \right| \\
&\le \frac{m\gamma}{\log \alpha} \sum_{i=f}^{f'} \frac{1}{i} \\
&\le \frac{C(1+\log f')}{\log \alpha}.
\end{aligned}
$$

To show (4.41), we calculate directly. For some $C$ (which may change from occurrence to occurrence) we obtain, noting that $\mathbf{P}[F > \lfloor \alpha \rfloor] = \mathcal{O}(\frac{1}{\alpha})$ from (4.17) and, for $i \le f$, $K_i$ and $L$ are conditionally independent given $F = f$,

$$
\begin{aligned}
&|\mathbf{P}[L=l|K_i=k] - \mathbf{P}[L=l]| \\
&= \left| \sum_{f=i}^{\lfloor \alpha \rfloor} (\mathbf{P}[L=l|F=f] \cdot \mathbf{P}[F=f|K_i=k] - \mathbf{P}[L=l] \cdot \mathbf{P}[F=f|K_i=k]) \right| \\
&\quad + \mathcal{O}\left( \frac{1}{\alpha} \right) \\
&= \left| \sum_{f=i}^{\lfloor \alpha \rfloor} \sum_{f'=n}^{\lfloor \alpha \rfloor} \mathbf{P}[F=f|K_i=k] \right. \\
&\qquad \left. \times \mathbf{P}[F=f'](\mathbf{P}[L=l|F=f] - \mathbf{P}[L=l|F=f']) \right| + \mathcal{O}\left( \frac{1}{\alpha} \right) \\
&\le C \sum_{f=i}^{\lfloor \alpha \rfloor} \sum_{f'=n}^{\lfloor \alpha \rfloor} \frac{i}{f^2} \frac{1}{f'^2} \frac{1(1+\log f)(1+\log f')}{\log \alpha} \\
&= \frac{C}{\log \alpha} \sum_{f=i}^{\lfloor \alpha \rfloor} \frac{i(1+\log f)}{f^2} = \frac{C(1+\log i)}{\log \alpha},
\end{aligned}
$$

where we have used again (4.17) and the integral (4.23).  $\square$

With the help of the previous lemma, we can now complete the proof of Proposition 3.11. First take $s > 0$. Then the assertion follows from the above statements, since by (3.5),

$$
\begin{aligned}
\mathbf{P}[L=l, S=s] \;=\;& \mathbf{P}[L=l, S=s, M=1] + \mathcal{O}\!\left(\frac{1}{(\log \alpha)^2}\right) \\[2mm]
=\;& \sum_{i=1}^{\lfloor \alpha \rfloor} \sum_{k=1}^{i \wedge (n-1)} \mathbf{P}[L=l \mid S=s, K_i=k, M=M_i=1] \\
& \qquad \times \mathbf{P}[S=s, K_i=k, M=M_i=1] + \mathcal{O}\!\left(\frac{1}{(\log \alpha)^2}\right) \\[2mm]
=\;& \sum_{i=1}^{\lfloor \alpha \rfloor} \sum_{k=1}^{i \wedge (n-1)} \mathbf{P}[L=l \mid K_i=k]\mathbf{P}[S=s, K_i=k, M=M_i=1] \\
& + \mathcal{O}\!\left(\frac{1}{(\log \alpha)^2}\right) \\[2mm]
\overset{(4.41)}{=}\;& \sum_{i=1}^{\lfloor \alpha \rfloor} \sum_{k=1}^{i \wedge (n-1)} \left(\mathbf{P}[L=l] + \frac{C \log i}{\log \alpha}\right) \\
& \qquad \times \mathbf{P}[S=s, K_i=k, M=M_i=1] + \mathcal{O}\!\left(\frac{1}{(\log \alpha)^2}\right) \\[2mm]
\overset{(4.21)}{=}\;& \mathbf{P}[L=l] \cdot \mathbf{P}[S=s] + \mathcal{O}\!\left(\frac{1}{(\log \alpha)^2}\right).
\end{aligned}
$$

Also for $s=0$ the assertion is true, since by the above

$$
\begin{aligned}
\mathbf{P}[L=l, S=0] &= \mathbf{P}[L=l] - \mathbf{P}[L=l, S>0] \\
&= \mathbf{P}[L=l](1 - \mathbf{P}[S>0]) + \mathcal{O}\!\left(\frac{1}{(\log \alpha)^2}\right) \\
&= \mathbf{P}[L=l] \cdot \mathbf{P}[S=0] + \mathcal{O}\!\left(\frac{1}{(\log \alpha)^2}\right). \qquad \square
\end{aligned}
$$

PROOF OF PROPOSITION 2.8.  Schweinsberg and Durrett [21] also used Yule processes to obtain an approximate scenario for the genealogy under hitchhiking. In fact, the third step in the approximation in our paper (see Sections 3.2 and 3.5) leads to the very same Yule process that appeared in [21]. To be exact about this, note that by Remark 3.7, one line in our Yule tree is hit by a mark during Yule time $i$ with probability

$$
\begin{aligned}
\frac{\rho}{i\alpha + \rho} &= \frac{\gamma}{\log \alpha} \frac{1}{i + \gamma/\log \alpha} \\
&= \frac{\gamma}{\log \alpha}\left(\frac{1}{i} + \mathcal{O}\!\left(\frac{1}{i^2 \log \alpha}\right)\right).
\end{aligned}
$$

In the model of Schweinsberg and Durrett [21] as given in their (7.2), a line is hit during Yule time $i$, as $\rho = 2Nr, \alpha = 2Ns$, with probability

$$\frac{r}{is + r(1-s)} = \frac{\rho}{i\alpha + \rho(1-s)}$$

$$= \frac{\gamma}{\log \alpha}\left(\frac{1}{i} + \mathcal{O}\left(\frac{1}{i^2 \log \alpha}\right)\right),$$

which proves the proposition. □

4.5. *The sampling formula.*

PROOF OF COROLLARY 2.7.   Using Theorem 1, we calculate, because $L$ and $S$ are independent,

$$\mathbf{P}[E = e, L = l] = \sum_{s=e}^{n} \mathbf{P}[E = e | L = l, S = s] \cdot \mathbf{P}[L = l, S = s]$$

$$= \mathbf{P}[L = l]\sum_{s=e}^{n} \frac{\binom{s}{e}\binom{n-s}{n-l-e}}{\binom{n}{l}}\mathbf{P}[S = s]$$

$$= \mathbf{E}[p_F^{n-l}(1-p_F)^l]\sum_{s=e}^{n}\binom{s}{e}\binom{n-s}{n-l-e}\mathbf{P}[S = s].$$

Let us now distinguish the three cases $e \geq 2, e = 1$ and $e = 0$. By a calculation using binomial coefficients (or using a computer algebra program or by [20], page 7, (2)),

$$\sum_{s=e}^{n}\binom{s-2}{e-2}\binom{n-s}{n-l-e} = \binom{n-1}{l}\mathbb{1}\{l + e \leq n\}.$$

With this we can calculate for $e \geq 2$, as long as $l + e \leq n$,

$$\sum_{s=e}^{n}\binom{s}{e}\binom{n-s}{n-l-e}\mathbf{P}[S = s]$$

$$= \frac{n\gamma}{\log \alpha}\left(\frac{1}{n}\binom{n}{e}\mathbb{1}\{l + e = n\} + \sum_{s=e}^{n}\binom{s}{e}\binom{n-s}{n-l-e}\frac{1}{s(s-1)}\right)$$

$$= \frac{n\gamma}{\log \alpha}\left(\frac{n-1}{e(e-1)}\binom{n-2}{e-2}\mathbb{1}\{l + e = n\}\right.$$

$$\left. + \frac{1}{e(e-1)}\sum_{s=e}^{n}\binom{s-2}{e-2}\binom{n-s}{n-l-e}\right)$$

$$= \frac{n\gamma}{\log \alpha}\frac{(n-1)\binom{n-2}{e-2}\mathbb{1}\{l + e = n\} + \binom{n-1}{l}}{e(e-1)},$$

where we have used $\frac{1}{n-1} = \frac{1}{n(n-1)} + \frac{1}{n}$ for the case $s = n$. This gives (2.7) in the case $e \geq 2$. For $e = 1$, we have

$$\sum_{s=1}^{n} \binom{s}{1} \binom{n-s}{n-l-1} \mathbf{P}[S = s]$$

$$= \frac{n\gamma}{\log \alpha} \left( \mathbb{1}\{l + 1 = n\} + \binom{n-1}{l} \sum_{i=2}^{n-1} \frac{1}{i} + \sum_{s=2}^{n} \frac{\binom{n-s}{l-s+1}}{s-1} \right),$$

which gives the result for $e = 1$. For the case $e = 0$, we first calculate

$$\mathbf{P}[S = 0] = 1 - \sum_{s=1}^{n} \mathbf{P}[S = s]$$

$$= 1 - \frac{n\gamma}{\log \alpha} \left( \frac{1}{n} + \sum_{i=2}^{n-1} \frac{1}{i} + \sum_{s=2}^{n} \frac{1}{s(s-1)} \right)$$

$$= 1 - \frac{n\gamma}{\log \alpha} \sum_{i=1}^{n-1} \frac{1}{i}.$$

With this we see

$$\sum_{s=0}^{n} \binom{n-s}{n-l} \mathbf{P}[S = s]$$

$$= \binom{n}{l} \left( 1 - \frac{n\gamma}{\log \alpha} \sum_{i=1}^{n-1} \frac{1}{i} \right)$$

$$+ \frac{n\gamma}{\log \alpha} \left( \frac{1}{n} \mathbb{1}\{l = n\} + \binom{n-1}{l-1} \sum_{i=2}^{n-1} \frac{1}{i} + \sum_{s=2}^{n} \binom{n-s}{n-l} \frac{1}{s(s-1)} \right)$$

$$= \binom{n}{l} \left( 1 - \frac{n\gamma}{\log \alpha} \left( 1 - \frac{l}{n} \sum_{i=2}^{n-1} \frac{1}{i} \right) \right)$$

$$+ \frac{n\gamma}{\log \alpha} \left( \frac{1}{n} \mathbb{1}\{l = n\} + \sum_{s=2}^{n} \binom{n-s}{n-l} \frac{1}{s(s-1)} \right),$$

which completes the proof. $\quad \square$

**Acknowledgments.** We thank Bob Griffiths, Jason Schweinsberg, Nick Barton and Wolfgang Angerer for illuminating and stimulating discussions. We gratefully acknowledge the hospitality of the Erwin Schrödinger Institute in Vienna, where part of our work was done. For hospitality during their research visit in February 2005, P. Pfaffelhuber and A. Wakolbinger thank the University of Oxford. Finally, we thank a referee and an Associate Editor for a very careful reading and helpful remarks.

## REFERENCES

[1] BARTON, N. (1998). The effect of hitch-hiking on neutral genealogies. *Gen. Res.* **72** 123–133.

[2] BARTON, N. H., ETHERIDGE, A. M. and STURM, A. (2004). Coalescence in a random background. *Ann. Appl. Probab.* **14** 754–785. MR2052901

[3] BIRKNER, M., BLATH, J., CAPALDO, A., ETHERIDGE, A., MÖHLE, M., SCHWEINSBERG, J. and WAKOLBINGER, A. (2005). Alpha-stable branching and beta coalescents. *Electron. J. Probab.* **10** 303–325. MR2120246

[4] DURRETT, R. and SCHWEINSBERG, J. (2004). Approximating selective sweeps. *Theor. Popul. Biol.* **66** 129–138.

[5] DURRETT, R. and SCHWEINSBERG, J. (2005). A coalescent model for the effect of advantageous mutations on the genealogy of a population. *Stochastic Process. Appl.* **115** 1628–1657. MR2165337

[6] ETHIER, S. N. and KURTZ, T. G. (1986). *Markov Processes*: *Characterization and Convergence.* Wiley, New York. MR0838085

[7] EVANS, S. N. and O'CONNELL, N. (1994). Weighted occupation time for branching particle systems and a representation for the supercritical superprocess. *Canad. Math. Bull.* **37** 187–196. MR1275703

[8] EWENS, W. J. (2004). *Mathematical Population Genetics. I. Theoretical Introduction*, 2nd ed. Springer, New York. MR2026891

[9] FISHER, R. A. (1930). *The Genetical Theory of Natural Selection*, 2nd ed. Clarendon Press.

[10] FRISTEDT, B. and GRAY, L. (1997). *A Modern Approach to Probability Theory.* Birkhäuser, Boston. MR1422917

[11] GRIFFITHS, R. C. (2003). The frequency spectrum of a mutation and its age, in a general diffusion model. *Theor. Popul. Biol.* **64** 241–251.

[12] GRIFFITHS, R. C. and TAVARÉ, S. (1998). The age of a mutation in a general coalescent tree. *Stochastic Models* **14** 273–295. MR1617552

[13] JOHNSON, N. L. and KOTZ, S. (1977). *Urn Models and Their Application.* Wiley, New York. MR0488211

[14] KAJ, I. and KRONE, S. M. (2003). The coalescent process in a population with stochastically varying size. *J. Appl. Probab.* **40** 33–48. MR1953766

[15] KAPLAN, N. L., HUDSON, R. R. and LANGLEY, C. H. (1989). The 'Hitchhiking effect' revisited. *Genetics* **123** 887–899.

[16] KARLIN, S. and TAYLOR, H. M. (1981). *A Second Course in Stochastic Processes.* Academic Press, New York. MR0611513

[17] KURTZ, T. G. (1971). Limit theorems for sequences of jupmp Markov processes approximating ordinary differential equations. *J. Appl. Probab.* **8** 344–356. MR0287609

[18] MAYNARD SMITH, J. and HAIGH, J. (1974). The hitch-hiking effect of a favorable gene. *Gen. Res.* **23** 23–35.

[19] O'CONNELL, N. (1993). Yule process approximation for the skeleton of a branching process. *J. Appl. Probab.* **30** 725–729. MR1232747

[20] RIORDAN, J. (1968). *Combinatorial Identities.* Wiley, New York. MR0231725

[21] SCHWEINSBERG, J. and DURRETT, R. (2005). Random partitions approximating the coalescence of lineages during a selective sweep. *Ann. Appl. Probab.* **15** 1591–1651. MR2152239

[22] STEPHAN, W., WIEHE, T. H. E. and LENZ, M. W. (1992). The effect of strongly selected substitutions on neutral polymorphism: Analytical results based on diffusion theory. *Theor. Popul. Biol.* **41** 237–254.

A. Etheridge
Department of Statistics
University of Oxford
1 South Parks Road
Oxford OX1 3TG
United Kingdom
E-mail: etheridg@stats.ox.ac.uk
URL: www.stats.ox.ac.uk/~etheridg

P. Pfaffelhuber
Department Biologie II
Ludwig-Maximilian University
Grosshaderner Strasse 2
82131 Planegg-Martinsried
Germany
E-mail: p.p@lmu.de
URL: www.zi.biologie.uni-muenchen.de/~pfaffelhuber

A. Wakolbinger
J. W. Goethe-Universität
Fachbereich Informatik und Mathematik
60054 Frankfurt am Main
Germany
E-mail: wakolbinger@math.uni-frankfurt.de
URL: ismi.math.uni-frankfurt.de/wakolbinger